\documentclass[a4paper,12pt,intlimits,oneside]{amsart}

\usepackage{enumerate}
\usepackage{amsfonts,amsmath}

\usepackage{latexsym,amssymb}

\newcommand{\comment}[1]{}

\newcounter{rea}
\newcounter{rek}
\newcounter{res}
\begin{document}

\title[New Hardy spaces of Musielak-Orlicz type]{New Hardy spaces of Musielak-Orlicz type and boundedness of sublinear operators}         

\author{Luong Dang  KY}    
\address{Department of Mathematics, University of Quy Nhon, 170 An Duong Vuong, Quy Nhon, Binh Dinh, Viet Nam} 
\email{{\tt dangky@math.cnrs.fr}}
\keywords{Muckenhoupt weights, Musielak-Orlicz functions, BMO-multipliers, Hardy spaces, atomic decompositions, Hardy-Orlicz spaces, quasi-Banach spaces, sublinear operators}
\subjclass[2010]{42B35 (46E30, 42B15, 42B30)}

\begin{abstract} We introduce a new class of Hardy spaces $H^{\varphi(\cdot,\cdot)}(\mathbb R^n)$,  called Hardy spaces of Musielak-Orlicz type, which generalize the Hardy-Orlicz spaces of Janson and the weighted Hardy spaces of Garc\'ia-Cuerva, Str\"omberg, and Torchinsky. Here, $\varphi: \mathbb R^n\times [0,\infty)\to [0,\infty)$ is a function such that $\varphi(x,\cdot)$ is an Orlicz function and $\varphi(\cdot,t)$ is a  Muckenhoupt $A_\infty$ weight. A function $f$ belongs to $H^{\varphi(\cdot,\cdot)}(\mathbb R^n)$ if and only if its maximal function $f^*$ is so that $x\mapsto \varphi(x,|f^*(x)|)$ is integrable.
Such a space arises naturally for instance in the description of the product of functions in $H^1(\mathbb R^n)$ and $BMO(\mathbb R^n)$ respectively (see \cite{BGK}).
 We characterize these spaces via the grand maximal function and establish their atomic decomposition. We characterize also their dual spaces.  The class of pointwise multipliers for $BMO(\mathbb R^n)$ characterized by Nakai and Yabuta can be seen as the dual of $L^1(\mathbb R^n)+  H^{\rm log}(\mathbb R^n)$ where $ H^{\rm log}(\mathbb R^n)$ is the Hardy space of Musielak-Orlicz type related to the Musielak-Orlicz function $\theta(x,t)=\displaystyle\frac{t}{\log(e+|x|)+ \log(e+t)}$.

Furthermore, under additional assumption on $\varphi(\cdot,\cdot)$ we prove that if $T$ is a sublinear operator and maps all atoms into uniformly bounded elements of a quasi-Banach space $\mathcal B$, then $T$ uniquely extends to a bounded sublinear operator from $H^{\varphi(\cdot,\cdot)}(\mathbb R^n)$ to $\mathcal B$. These results are new even for the classical  Hardy-Orlicz spaces on $\mathbb R^n$.
\end{abstract}

\maketitle
\newtheorem{theorem}{Theorem}[section]
\newtheorem{lemma}{Lemma}[section]
\newtheorem{proposition}{Proposition}[section]
\newtheorem{remark}{Remark}[section]
\newtheorem{corollary}{Corollary}[section]
\newtheorem{definition}{Definition}[section]
\newtheorem{example}{Example}[section]
\numberwithin{equation}{section}
\newtheorem{Theorem}{Theorem}[section]
\newtheorem{Lemma}{Lemma}[section]
\newtheorem{Proposition}{Proposition}[section]
\newtheorem{Remark}{Remark}[section]
\newtheorem{Corollary}{Corollary}[section]
\newtheorem{Definition}{Definition}[section]
\newtheorem{Example}{Example}[section]
\newtheorem*{lemmaa}{Lemma A}
\newtheorem*{lemmab}{Lemma B}
\newtheorem*{lemmac}{Lemma C}
\newtheorem*{lemmad}{Lemma D}
\newtheorem*{lemmae}{Lemma E}
\newtheorem*{lemmaf}{Lemma F}

\section{Introduction}

Since Lebesgue theory of integration has taken a center stage in concrete problems of analysis, the need for more inclusive classes of function spaces than the $L^p(\mathbb R^n)$-families naturally arose. It is well known that the Hardy spaces $H^p(\mathbb R^n)$ when $p\in (0,1]$ are good substitutes of $L^p(\mathbb R^n)$ when studying the boundedness of operators: for example, the Riesz operators are  bounded on $H^p(\mathbb R^n)$, but not on $L^p(\mathbb R^n)$ when $p\in (0,1]$. The theory of Hardy spaces $H^p$ on the Euclidean space $\mathbb R^n$ was initially developed by Stein and Weiss \cite{SW}. Later, Fefferman and Stein \cite{FeS} systematically developed a real-variable theory for the Hardy spaces $H^p(\mathbb R^n)$ with $p\in (0,1]$, which now plays an important role in various fields of analysis and partial differential equations; see, for example, \cite{CW, CLMS, Mul}. A key feature of the classical Hardy spaces is their atomic decomposition characterizations, which were obtained by Coifman \cite{Co} when $n = 1$ and Latter \cite{La} when $n > 1$. Later, the theory of Hardy spaces and their dual spaces associated with Muckenhoupt weights have been extensively studied by Garc\'ia-Cuerva \cite{Ga}, Str\"omberg and Torchinsky \cite{ST} (see also \cite{MW2, Bu, GM}); there the weighted Hardy spaces was defined by using the nontangential maximal functions and the atomic decompositions were derived. On the other hand, as another generalization of $L^p(\mathbb R^n)$, the Orlicz spaces were introduced by Birnbaum-Orlicz in \cite{BO} and Orlicz in \cite{Or}, since then, the theory of the Orlicz spaces themselves has been well developed and the spaces have been widely used in probability, statistics, potential theory, partial differential equations, as well as harmonic analysis and some other fields of analysis; see, for example, \cite{AIKM, IO, MaW}. Moreover, the Hardy-Orlicz spaces are also good substitutes of the Orlicz spaces in dealing with many problems of analysis, say, the boundedness of operators. 

Let $\Phi$ be a Orlicz function which is of positive lower type and (quasi-)concave. In \cite{Ja2}, Janson has considered the Hardy-Orlicz space $H^\Phi(\mathbb R^n)$   the space of all tempered distributions $f$ such that the nontangential grand maximal function of $f$ is defined by
$$f^*(x)=\sup\limits_{\phi\in \mathcal A_N}\sup\limits_{|x-y|<t}|f*\phi_t(y)|,$$
for all $x\in\mathbb R^n$, here  and in what follows $\phi_t(x):= t^{-n}\phi(t^{-1}x)$, with
$$\mathcal A_N= \Big\{\phi\in \mathcal S(\mathbb R^n): \sup\limits_{x\in\mathbb R^n}(1+ |x|)^N |\partial^\alpha_x \phi(x)|\leq 1\; \mbox{for}\; \alpha\in \mathbb N^n, |\alpha|\leq N\Big\}$$
with $N=N(n,\Phi)$   taken large enough, belongs to the Orlicz space $L^\Phi(\mathbb R^n)$.  Recently, the theory of  Hardy-Orlicz spaces associated with operators (see \cite{CCYY1, CCYY2, JY, YY3}) have also been introduced and studied. Remark that these Hardy-Orlicz type spaces appear naturally when studying  the theory of nonlinear PDEs  (cf. \cite{GI, IS, IV2}) since many cancellation phenomena for Jacobians cannot be observed in the usual Hardy spaces $H^p(\mathbb R^n)$. For instance,  let $f=(f^1,..., f^n)$ in the Sobolev class $W^{1,n}(\mathbb R^n, \mathbb R^n)$ and the Jacobians $J(x,f)dx= df^1\wedge\cdots\wedge df^n$, then (see  Theorem 10.2 of \cite{IV2})
$$\mathcal T(J(x,f))\in L^1(\mathbb R^n)+ H^\Phi(\mathbb R^n)$$
where $\Phi(t)=t/\log(e+t)$ and $\mathcal T(f)= f\log|f|$, since $J(x,f)\in H^1(\mathbb R^n)$ (cf. \cite{CLMS}) and $\mathcal T$ is well defined on $H^1(\mathbb R^n)$. We refer  readers to \cite{RW, IV1}  for this interesting nonlinear operator $\mathcal T$.

In this paper we want to allow generalized Hardy-Orlicz spaces related to generalized Orlicz functions that may vary in the spatial variables. More precisely the Orlicz function $\Phi(t)$ is replaced by a function $\varphi(x,t)$, called Musielak-Orlicz function (cf. \cite{Mus, Di}). We then define  Hardy spaces of Musielak-Orlicz type. Apart from interesting theoretical considerations, the motivation to study function spaces of Musielak-Orlicz type comes from applications to elasticity, fluid dynamics, image processing, nonlinear PDEs and the calculus of variation (cf. \cite{Di, DHR}).

A particular case of Hardy spaces of Musielak-Orlicz type appears naturally when considering the  products of   functions in $BMO(\mathbb R^n)$ and   $ H^1(\mathbb R^n)$ (see \cite{BGK}); and the endpoint estimates for the div-curl lemma (see \cite{BFG, BGK}). More precisely, in \cite{BGK} the authors proved that product of a  $BMO(\mathbb R^n)$ function and a $ H^1(\mathbb R^n)$ function may be written as a sum of an integrable term and of a term in $H^{\log}(\mathbb R^n)$,  a Hardy space of Musielak-Orlicz type related to the Musielak-Orlicz function $\theta(x,t)=\frac{t}{\log(e+|x|)+ \log(e+t)}$. Moreover, the corresponding bilinear operators are bounded. This result  gives in particular a positive answer to the Conjecture 1.7  in \cite{BIJZ}. By duality, one finds  pointwise multipliers for $BMO(\mathbb R^n)$. Recall that a function $g$ on $\mathbb R^n$ is called a pointwise multiplier for $BMO(\mathbb R^n)$, if the pointwise multiplication $fg$ belongs to $BMO(\mathbb R^n)$ for all $f$ in  $BMO(\mathbb R^n)$.  In \cite{NY}, Nakai and Yabuta characterize the pointwise multipliers for $BMO(\mathbb R^n)$: they prove that $g$ is a  pointwise multiplier for $BMO(\mathbb R^n)$ if and only if $g$ belong to $L^\infty(\mathbb R^n)\cap BMO^{\rm log}(\mathbb R^n)$, where $BMO^{\rm log}(\mathbb R^n)$ is the space of all locally integrable functions $f$ such that
\begin{eqnarray*}
 \|f\|_{BMO^{\rm log}}:= \sup\limits_{B(a,r)}\frac{|\log r|+ \log(e+|a|)}{|B(a,r)|}\int_{B(a,r)}|f(x)-f_{B(a,r)}|dx<\infty.
\end{eqnarray*}
By  using the theory of these new Hardy spaces and dual spaces, we establish that the class of pointwise multipliers for $BMO(\mathbb R^n)$   is just the dual of $L^1(\mathbb R^n)+  H^{\rm log}(\mathbb R^n)$. Remark that the class of pointwise multipliers for $BMO(\mathbb R^n)$  have also recently been used by Lerner \cite{Le} for solving a  conjecture of  Diening (see \cite{Di}) on the boundedness of the Hardy-Littlewood maximal operator on the generalized Lebesgue spaces $L^{p(x)}(\mathbb R^n)$ (a special case of Musielak-Orlicz spaces, for the details see \cite{Di, Le}).

Motivated by all of the above mentioned facts, in this paper, we introduce a new class of Hardy spaces $H^{\varphi(\cdot,\cdot)}(\mathbb R^n)$,  called Hardy spaces of Musielak-Orlicz type, which generalize the Hardy-Orlicz spaces of Janson and the weighted Hardy spaces of Garc\'ia-Cuerva, Str\"omberg, and Torchinsky. Here, $\varphi: \mathbb R^n\times [0,\infty)\to [0,\infty)$ is a function such that $\varphi(x,\cdot)$ is an Orlicz function and $\varphi(\cdot,t)$ is a  Muckenhoupt weight  $A_\infty$. In the special case  $\varphi(x,t)=w(x)\Phi(t)$ with $w$ in the Muckenhoupt class  and $\Phi$ an Orlicz function, our Hardy spaces are weighted Hardy-Orlicz spaces but they are different from the ones considered  by  Harboure, Salinas, and Viviani \cite{HSV1, HSV2}.

As an example of our results, let us give the atomic decomposition with bounded atoms. Let $\varphi$ be a {\sl growth function} (see Section 2). A bounded function $a$ is a $\varphi$-atom if it satisfies the following three conditions

i) supp $a\subset B$ for some ball $B$,

ii) $\|a\|_{L^\infty}\leq \|\chi_B\|_{L^\varphi}^{-1}$,

iii) $\int_{\mathbb R^n}a(x)x^\alpha dx=0$ for any $|\alpha|\leq [n(\frac{q(\varphi)}{i(\varphi)}-1)]$,\\
where $q(\varphi)$ and $i(\varphi)$ are the indices of $\varphi$ (see Section 2). We next define the {\sl atomic Hardy space of Musielak-Orlicz type} $H^{\varphi(\cdot,\cdot)}_{\rm at}(\mathbb R^n)$ as those distributions $f\in \mathcal S'(\mathbb R^n)$ such that $f=\sum_j  b_j$ (in the sense of $\mathcal S'(\mathbb R^n)$), where ${b_j}^,s$ are multiples of $\varphi$-atoms supported in the balls ${B_j}^,s$, with the property $\sum_j \varphi(B_j, \|b_j\|_{L^q_\varphi(B_j)})<\infty$; and define the norm of $f$ by
$$\|f\|_{ H^{\varphi(\cdot,\cdot)}_{\rm at}}=\inf\Big\{\Lambda_\infty(\{ b_j\}): f= \sum_j  b_j \quad\mbox{in the sense of}\;\;\mathcal S'(\mathbb R^n)\Big\},$$
where $\Lambda_\infty(\{ b_j\})=\inf\Big\{\lambda>0:\sum_j \varphi\Big(B_j, \frac{\|b_j\|_{L^\infty}}{\lambda}\Big)\leq 1\Big\}$ with $\varphi(B,t):= \int_B \varphi(x,t)dx$ for all $t\geq 0$ and $B$ is measurable. Then we obtain:
\begin{Theorem}
$H^{\varphi(\cdot,\cdot)}_{\rm at}(\mathbb R^n)= H^{\varphi(\cdot,\cdot)}(\mathbb R^n)$ with equivalent norms.
\end{Theorem}

The fact that $\Lambda_\infty(\{ b_j\})$, which is the right expression for the (quasi-)norm in the atomic Hardy space of Musielak-Orlicz type, plays a central role in this paper. It should be emphasized that, even if the steps of the proof of such a theorem are standard, the adaptation to this context is not standard.

On the other hand, to establish the boundedness of operators on Hardy spaces, one usually appeals to the atomic decomposition characterization, see \cite{Co, La, TW}, which means that a function or distribution in Hardy spaces can be represented as a linear combination of functions of an elementary form, namely, atoms. Then, the boundedness of operators on Hardy spaces can be deduced from their behavior on atoms or molecules in principle. However, caution needs to be taken due to an example constructed in Theorem 2 of \cite{Bo}. There exists a linear functional defined on a dense subspace of $H^1(\mathbb R^n)$, which maps all $(1,\infty,0)$-atoms into bounded scalars, but however does not extend to a bounded linear functional on the whole $H^1(\mathbb R^n)$. This implies that the uniform boundedness of a linear operator $T$ on atoms does not automatically guarantee the boundedness of $T$ from $H^1(\mathbb R^n)$ to a Banach space $\mathcal B$. Nevertheless, by using the grand maximal function characterization of $H^p(\mathbb R^n)$, Meda, Sj\"ogren, and Vallarino \cite{MSV1, MSV2} proved that if a sublinear operator $T$ maps all $(p,q,s)$-atoms when $q<\infty$ and continuous $(p,\infty,s)$-atoms into uniformly bounded elements of $L^p(\mathbb R^n)$ (see also \cite{YZ, BLYZ} for quasi-Banach spaces), then $T$ uniquely extends to a bounded sublinear operator from $H^p(\mathbb R^n)$ to $L^p(\mathbb R^n)$. In this paper, we study boundedness of sublinear operators in the context of new Hardy spaces of Musielak-Orlicz type which generalize the main results in \cite{MSV1, MSV2}. More precisely, under additional assumption on $\varphi(\cdot,\cdot)$, we prove that finite atomic norms on dense subspaces of  $H^{\varphi(\cdot,\cdot)}(\mathbb R^n)$ are equivalent with the standard infinite atomic decomposition norms. As an application, we prove that if $T$ is a sublinear operator and maps all atoms into uniformly bounded elements of a quasi-Banach space $\mathcal B$, then $T$ uniquely extends to a bounded sublinear operator from $H^{\varphi(\cdot,\cdot)}(\mathbb R^n)$ to $\mathcal B$.

Using the theory of these new Hardy spaces and ideas from \cite{BGK},  we studied and established (see \cite{Ky2, Ky3}) some new interesting estimates of endpoint type for the commutators of  singular integrals and fractional integrals on Hardy-type spaces. Recently, new Hardy spaces of Musielak-Orlicz type have developed in many directions (see \cite{BCKYY, NS, Sa, YY, YY2}), some their charcterizations and applications were also established and studied in \cite{BCKYY, Sa, HYY, LHY}.

Our paper is organized as follows. In Section 2 we give the notation and definitions that we shall use in the sequel. For simplicity we  write $\varphi$ for $\varphi(\cdot, \cdot)$. One then introduces Hardy spaces of Musielak-Orlicz type $H^\varphi(\mathbb R^n)$ via grand maximal functions, atomic Hardy spaces $H^{\varphi, q, s}_{\rm at}(\mathbb R^n)$, finite  atomic Hardy spaces $H^{\varphi, q, s}_{\rm fin}(\mathbb R^n)$ for any admissible triplet $(\varphi,q,s)$, $BMO$-Musielak-Orlicz-type spaces  $BMO^{\varphi}(\mathbb R^n)$, and generalized quasi-Banach spaces $\mathcal B_\gamma$ for $\gamma\in (0,1]$. In Section 3 we state the main results: the atomic decompositions (Theorem \ref{atomic decomposition for new Hardy spaces}), the duality (Theorem \ref{the dual theorem for new Hardy spaces}), the class of pointwise multipliers for $BMO(\mathbb R^n)$ (Theorem \ref{a theorem on the class of pointwise multipliers for BMO}), the finite atomic decomposition (Theorem \ref{finite decomposition for new Hardy spaces}), and the criterion for boundedness of sublinear operators in $H^\varphi(\mathbb R^n)$ (Theorem \ref{the boundedness of sublinear operators on new Hardy spaces}). In Section 4 we present and prove the basic properties of the growth functions $\varphi$ since they provide the tools for further work with this type of functions. In Section 5 we generalize the Calder\'on-Zygmund decomposition associated to the grand maximal function on $\mathbb R^n$ in the setting of the spaces of Musielak-Orlicz type. Applying this, we further prove that for any admissible triplet $(\varphi, q,s)$, $H^\varphi(\mathbb R^n)= H^{\varphi,q,s}_{\rm at}(\mathbb R^n)$ with equivalent norms (Theorem \ref{atomic decomposition for new Hardy spaces}). In Section 6 we prove the dual theorem. By Theorem 2 in \cite{Bo}, one has to be careful with the argument {\sl "the operator $T$ is uniformly bounded in $H^p_w(\mathbb R^n)$ ($H^\varphi(\mathbb R^n)$ here $\varphi(x,t)=w(x).t^p$ in our context) on $w$-$(p,\infty)$-atoms, and hence it extends to a bounded operator on $H^p_w(\mathbb R^n)$"} which has been used in \cite{Ga} and \cite{Bu}. In Section 7 we introduce {\sl log-atoms} and consider the particular case of $H^{\rm log}(\mathbb R^n)$. Finally, in Section 8 we  prove that $\|\cdot\|_{H^{\varphi,q,s}_{\rm fin}}$ and $\|\cdot\|_{H^\varphi}$ are equivalent quasi-norms on $H^{\varphi,q,s}_{\rm fin}(\mathbb R^n)$ when $q<\infty$ and on $H^{\varphi,q,s}_{\rm fin}(\mathbb R^n)\cap C(\mathbb R^n)$ when $q=\infty$, here and in  what follows  $C(\mathbb R^n)$ denotes the set of all continuous functions. Then,  we consider generalized quasi-Banach spaces which generalize the notion of quasi-Banach spaces in \cite{YZ} (see also \cite{BLYZ}), and obtain criterious for boundedness of sublinear operators on $H^\varphi(\mathbb R^n)$.

Throughout the whole paper, $C$ denotes a positive geometric constant which is independent of the main parameters, but may change from line to line.  The symbol $f\approx g$ means that $f$ is equivalent to $g$ (i.e. $C^{-1}f\leq g\leq C f$), and $[\cdot]$ denotes the integer function. By $X^*$ we denote the dual of the (quasi-)Banach space $X$. In $\mathbb R^n$, we denote by $B=B(x,r)$ an open ball with center $x$ and radius $r>0$. For any measurable set $E$, we denote by $\chi_E$ its characteristic function, by $|E|$ its  Lebesgue measure, and by $E^c$ the set $\mathbb R^n\setminus E$.

\section{Notation and definitions}\label{the section of notation and definitions}

\subsection{Musielak-Orlicz-type functions}

First let us recall notations for Orlicz functions.

A function $\phi:[0,\infty)\to [0,\infty)$ is called {\sl Orlicz} if it is nondecreasing and $\phi(0)=0$; $\phi(t)>0, t>0$;  $\lim_{t\to \infty}\phi(t)=\infty$. An Orlicz function $\phi$ is said to be of {\sl lower type} (resp., {\sl upper type}) $p$, $p\in (-\infty,\infty)$, if there exists a positive constant $C$ so that
$$\phi(st)\leq C s^p\phi(t),$$
for all $t\geq 0$ and $s\in (0,1)$ (resp., $s\in [1,\infty)$). One say that $\phi$ is of {\sl positive lower type} (resp., {\sl finite upper type}) if it is of lower type (resp., upper type) $p$ for some $p>0$ (resp., $p$ finite).

Obviously, if $\phi$ is both of lower type $p_1$ and of upper type $p_2$, then $p_1\leq p_2$. Moreover, if $\phi$ is of lower type (resp., upper type) $p$ then it is also of lower type (resp., upper) $\widetilde p$ for $-\infty< \widetilde p<p$ (resp., $p<\widetilde p<\infty$). We thus write $$i(\phi):=\sup\{p\in (-\infty,\infty): \phi\;\mbox{is of lower type}\; p\}$$
$$I(\phi):= \inf\{p\in (-\infty,\infty): \phi\;\mbox{is of upper type}\; p\}$$
to denote the critical lower type and the critical upper type of the function $\phi$. 

Let us generalize these notions to functions $\varphi:\mathbb R^n\times [0,\infty)\to [0,\infty)$.

Given a function $\varphi:\mathbb R^n\times [0,\infty)\to [0,\infty)$ so that for any $x\in \mathbb R^n$, $\varphi(x,\cdot)$ is Orlicz. We say that  $\varphi$ is of {\sl uniformly lower type} (resp., {\sl upper type}) $p$ if there exists a positive constant $C$ so that
\begin{equation}\label{lower type}
\varphi(x, st)\leq C s^p\varphi(x, t),
\end{equation}
for all $x\in \mathbb R^n$ and $t\geq 0, s\in (0,1)$ (resp., $s\in [1,\infty)$). We say that $\varphi$ is of {\sl positive  uniformly lower type} (resp., {\sl finite uniform upper type}) if it is of uniformly lower type (resp., uniform upper type) $p$ for some $p>0$ (resp., $p$ finite), and denote
$$i(\varphi):= \sup\{p\in (-\infty,\infty): \varphi\;\mbox{is of uniformly lower type}\;p\}$$
$$I(\varphi):= \inf\{p\in (-\infty,\infty): \varphi\;\mbox{is of uniformly upper type}\;p\}.$$

We next  need to recall  notations for Muckenhoupt weights.

Let $1\leq q<\infty$. A nonnegative locally integrable function $w$ belongs to the {\sl Muckenhoupt class} $A_q$, say $w\in A_q$, if there exists a positive constant $C$ so that
$$\frac{1}{|B|}\int_B w(x)dx\Big(\frac{1}{|B|}\int_B (w(x))^{-1/(q-1)}dx\Big)^{q-1}\leq C, \quad\mbox{if}\; 1<q<\infty,$$
and 
$$\frac{1}{|B|}\int_B w(x)dx\leq C \mathop{\mbox{ess-inf}}\limits_{x\in B}w(x),\quad\mbox{if}\; q=1,$$
for all balls $B$ in $\mathbb R^n$. We say that $w\in A_\infty$ if $w\in A_q$ for some $q\in [1,\infty)$.

It is well known that $w\in A_q$, $1\leq q<\infty$, implies $w\in A_r$ for all $r >q$. Also, if $w\in A_q$, $1<q<\infty$, then $w\in A_r$ for some $r\in [1,q)$. One thus write $q_w := \inf\{q \geq 1: w\in A_q \}$ to denote the critical index of $w$.

Now, let us generalize these notions to functions $\varphi:\mathbb R^n\times [0,\infty)\to [0,\infty)$.

Let $\varphi:\mathbb R^n\times [0,\infty)\to \mathbb C$ be so that $x\mapsto \varphi(x,t)$ is measurable for all $t\in [0,\infty)$. We say that  $\varphi(\cdot,t)$ is {\sl uniformly locally integrable}  if for all compact set $K$ in $\mathbb R^n$, the following holds
$$\int\limits_K \sup\limits_{t>0}\frac{|\varphi(x,t)|}{\int_K |\varphi(y,t)|dy}dx<\infty$$
whenever the integral  exists.
A simple example for such {\sl uniformly locally integrable} functions is $\varphi(x,t)=w(x)\Phi(t)$ with $w$  a locally integrable function on $\mathbb R^n$ and $\Phi$ an arbitrary function on $[0,\infty)$. Our interesting examples are {\sl uniformly locally integrable} functions $\varphi(x,t)=\frac{t^p}{(\log(e+|x|)+ \log(e+t^p))^p}, 0<p\leq 1$, since they arise naturally in the study of pointwise product of functions in $H^p(\mathbb R^n)$ with functions in $BMO(\mathbb R^n)$ (cf. \cite{BGK}).

Given $\varphi:\mathbb R^n\times [0,\infty)\to [0,\infty)$ is a {\sl uniformly locally integrable} function.  We say that the function $\varphi(\cdot,t)$ satisfies the {\sl uniformly Muckenhoupt} condition $\mathbb A_q$, say $\varphi\in\mathbb A_q$, for some $1\leq q< \infty$  if there exists a positive constant $C$ so that
$$\frac{1}{|B|}\int_B \varphi(x,t)dx. \Big(\frac{1}{|B|}\int_B \varphi(x,t)^{-1/(q-1)}dx\Big)^{q-1}\leq C, \quad\mbox{if}\; 1<q<\infty,$$
and
$$\frac{1}{|B|}\int_B \varphi(x,t)dx\leq C \mathop{\mbox{ess-inf}}\limits_{x\in B} \varphi(x,t),\quad\mbox{if}\; q=1,$$
for all $t>0$ and balls $B$ in $\mathbb R^n$. We also say that $\varphi\in \mathbb A_\infty$ if $\varphi\in \mathbb A_q$ for some $q\in [1, \infty)$, and denote 
$$q(\varphi):= \inf\{q\geq 1: \varphi\in \mathbb A_q\}.$$

Now, we are able to introduce the {\sl growth functions} which are the basis for our new Hardy spaces.
\begin{Definition}
We say that $\varphi:\mathbb R^n\times [0,\infty)\to [0,\infty)$ is a {\bf growth function} if the following conditions are satisfied.
\begin{enumerate}

\item  The function $\varphi$ is a {\bf Musielak-Orlicz function} that is
\begin{enumerate} 
\item the function $\varphi(x,\cdot): [0,\infty)\to [0,\infty)$  is an Orlicz function for all $x\in \mathbb R^n$,
\item the function $\varphi(\cdot,t)$ is a Lebesgue measurable function for all $t\in [0,\infty)$.
\end{enumerate}
\item The function $\varphi$ belongs to  $\mathbb A_\infty$.
\item The function $\varphi$ is of positive uniformly lower type and of uniformly upper type 1.
\end{enumerate}

For $\varphi$  a {\sl growth function}, we denote $m(\varphi):= \Big[n\Big(\frac{q(\varphi)}{i(\varphi)}-1\Big)\Big]$.
\end{Definition}

Clearly, $\varphi(x,t)=w(x)\Phi(t)$ is a {\sl growth function} if $w\in A_\infty$ and $\Phi$ is of positive lower type and of upper type 1. Of course, there exists {\sl growth functions} which are not of that form for instance $\varphi(x,t)=\frac{t^\alpha}{[\log(e+|x|)]^\beta+ [\log(e+t)]^\gamma}$ for $\alpha\in (0,1]; \beta, \gamma\in (0,\infty)$. More precisely,  $\varphi\in \mathbb A_1$ and $\varphi$ is of uniformly upper type $\alpha$ with $i(\varphi)=\alpha$. In this paper, we are especially interested in the {\sl growth functions} $\varphi(x,t)=\frac{t^p}{(\log(e+|x|)+ \log(e+t^p))^p}, 0<p\leq 1$, since the Hardy spaces of Musielak-Orlicz type $H^\varphi(\mathbb R^n)$ arise naturally in the study of pointwise product of functions in $H^p(\mathbb R^n)$ with functions in $BMO(\mathbb R^n)$ (see also \cite{BG} in the setting of holomorphic functions in convex domains of finite type or  strictly pseudoconvex domains in $\mathbb C^n$).

\subsection{Hardy spaces of Musielak-Orlicz type}

Throughout the whole paper, we always assume that $\varphi$ is a {\sl growth function}.

Let us now introduce the {\sl Musielak-Orlicz-type spaces}. 

The Musielak-Orlicz-type space $L^\varphi(\mathbb R^n)$ is the set of all measurable functions $f$ such that $\int_{\mathbb R^n}\varphi(x,|f(x)|/\lambda)dx<\infty$ for some $\lambda>0$, with Luxembourg (quasi-)norm
$$\|f\|_{L^\varphi}:= \inf\Big\{\lambda>0: \int_{\mathbb R^n}\varphi(x,|f(x)|/\lambda)dx\leq 1\Big\}.$$

As usual, $\mathcal S(\mathbb R^n)$ denote the Schwartz class of test functions on $\mathbb R^n$ and $\mathcal S'(\mathbb R^n)$ the space of tempered distributions (or distributions for brevity). For $m\in\mathbb N$, we define 
$$\mathcal S_m(\mathbb R^n)=\Big\{\phi\in\mathcal S(\mathbb R^n):\|\phi\|_m=\sup\limits_{x\in\mathbb R^n,|\alpha|\leq m+1}(1+|x|)^{(m+2)(n+1)}|\partial^\alpha_x\phi(x)|\leq 1\Big\}.$$

For each distribution $f$, we define the nontangential grand maximal function $f^*_m$ of $f$ by

$$f^*_m(x)=\sup\limits_{\phi\in \mathcal S_m(\mathbb R^n)}\sup\limits_{|y-x|<t}|f*\phi_t(y)|,\; x\in\mathbb R^n.$$

When $m=m(\varphi)$ we write $f^*$ instead of $f^*_{m(\varphi)}$.

\begin{Definition}
The Hardy space of Musielak-Orlicz type $H^\varphi(\mathbb R^n)$ is the space of all distributions $f$ such that $f^*\in L^\varphi(\mathbb R^n)$ with the (quasi-)norm
$$\|f\|_{H^\varphi}:= \|f^*\|_{L^\varphi}.$$
\end{Definition}

Observe that, when $\varphi(x,t)=w(x)\Phi(t)$ with $w$ a Muckenhoupt weight and $\Phi$  an Orlicz function, our Hardy spaces are weighted Hardy-Orlicz spaces which include the classical Hardy-Orlicz spaces of Janson \cite{Ja2} ($w\equiv 1$ in this context) and the classical weighted Hardy spaces of Garc\'ia-Cuerva \cite{Ga}, Str\"omberg and Torchinsky \cite{ST} ($\Phi(t)\equiv t^p$ in this context), see also \cite{MW2, Bu, GM} .

Next, to introduce the {\sl atomic Hardy spaces of Musielak-Orlicz type} below, we need the following new spaces.
\begin{Definition}
For each ball $B$ in $\mathbb R^n$, we denote $L^q_\varphi(B), 1\leq q\leq \infty$, the set of all measurable functions $f$ on $\mathbb R^n$ supported in $B$ such that
\begin{equation}
\|f\|_{L^q_\varphi(B)}:= 
\begin{cases}
\sup\limits_{t>0}\Big(\frac{\int_{\mathbb R^n}|f(x)|^q\varphi(x,t)dx}{\varphi(B,t)}\Big)^{1/q}<\infty &,\quad 1\leq q<\infty, \\
\|f\|_{L^\infty}<\infty &, \quad q=\infty,
\end{cases}
\end{equation}
here and in the future $\varphi(B,t):= \int_B \varphi(x,t)dx$.
\end{Definition}
Then, it is straightforward to verify that $(L^q_\varphi(B),\|\cdot\|_{L^q_\varphi(B)})$ is a Banach space.

Now, we are able to  introduce the {\sl atomic Hardy spaces of Musielak-Orlicz type}.
\begin{Definition}\label{atom}
A triplet $(\varphi,q,s)$ is called admissible, if $q\in(q(\varphi),\infty]$ and $s\in\mathbb N$ satisfies $s\geq m(\varphi)$. A measurable function $a$  is a $(\varphi,q,s)$-atom if it satisfies the following three conditions

i) $a\in L^q_\varphi(B)$ for some ball B,

ii) $\|a\|_{L^q_\varphi(B)}\leq  \|\chi_B\|^{-1}_{L^\varphi}$, 

iii) $\int_{\mathbb R^n}a(x)x^{\alpha}dx=0$ for any $|\alpha|\leq s$.
\end{Definition}

In this setting we define the {\sl atomic Hardy space of Musielak-Orlicz type} $H^{\varphi,q,s}_{\rm at}(\mathbb R^n)$ as those distributions $f\in \mathcal S'(\mathbb R^n)$ that can be represented as a sum of multiples of $(\varphi,q,s)$-atoms, that is,
$$f=\sum_j  b_j\quad\mbox{in the sense of}\;\;\mathcal S'(\mathbb R^n),$$
where ${b_j}^,s$ are multiples of $(\varphi,q,s)$-atoms supported in the balls ${B_j}^,s$, with the property
$$\sum_j \varphi(B_j, \|b_j\|_{L^q_\varphi(B_j)})<\infty.$$
We introduce a (quasi-)norm in $ H^{\varphi,q,s}_{\rm at}(\mathbb R^n)$. Given a sequence of multiples of $(\varphi,q,s)$-atoms, $\{b_j\}_j$, we denote
\begin{equation}
\Lambda_q(\{ b_j\})=\inf\Big\{\lambda>0:\sum_j \varphi\Big(B_j, \frac{\|b_j\|_{L^q_\varphi(B_j)}}{\lambda}\Big)\leq 1\Big\}
\end{equation}
and define
\begin{equation}
\|f\|_{ H^{\varphi,q,s}_{\rm at}}=\inf\Big\{\Lambda_q(\{ b_j\}): f= \sum_j  b_j \quad\mbox{in the sense of}\;\;\mathcal S'(\mathbb R^n)\Big\}.
\end{equation}

Let $(\varphi,q,s)$ be an admissible triplet. We denote $H^{\varphi,q,s}_{\rm fin}(\mathbb R^n)$ the vector space of all finite linear combinations of $(\varphi,q,s)$-atoms, that is,
$$ f=\sum_{j=1}^k b_j,$$
where $b_j$'s are multiples of $(\varphi,q,s)$-atoms supported in balls $B_j$'s. Then, the norm of $f$ in $H^{\varphi,q,s}_{\rm fin}(\mathbb R^n)$ is defined by
\begin{equation}
\|f\|_{H^{\varphi,q,s}_{\rm fin}}=\inf\Big\{\Lambda_q(\{b_j\}_{j=1}^k): f= \sum_{j=1}^k b_j\Big\}.
\end{equation}

Obviously, for any admissible triplet $(\varphi,q,s)$, the set $H^{\varphi,q,s}_{\rm fin}(\mathbb R^n)$ is dense in $H^{\varphi,q,s}_{\rm at}(\mathbb R^n)$ with respect to the quasi-norm $\|\cdot\|_{H^{\varphi,q,s}_{\rm at}}$.

We should point out that the theory of atomic Hardy-Orlicz spaces  have been first introduced by Viviani \cite{Vi} in the setting of spaces of homogeneous type. Later,  Serra \cite{Ser} generalized it to the context of the Euclidean space $\mathbb R^n$ and obtained the molecular characterization. In the particular case, when $\varphi(x,t)\equiv \Phi(t)$ the space $H^{\varphi,q,s}_{\rm at}(\mathbb R^n)$ is the space considered in \cite{Ser}. We also remark that when $\varphi(x,t)\equiv w(x).t^p, 0<p\leq 1$, $w$ a Muckenhoupt weight, the space $H^{\varphi,q,s}_{\rm at}(\mathbb R^n)$ is just the classical weighted atomic Hardy space $H^{p,q,s}_w(\mathbb R^n)$ which has been considered by Garc\'ia-Cuerva \cite{Ga}, Str\"omberg and Torchinsky \cite{ST}.

\subsection{BMO-Musielak-Orlicz-type spaces}

We also need $BMO$ type spaces, which will be in duality of the Hardy spaces of Musielak-Orlicz type  defined above.  A function $f\in L^1_{\rm loc}(\mathbb R^n)$ is said to belong to $BMO^{\varphi}(\mathbb R^n)$ if
$$\|f\|_{BMO^\varphi}:= \sup\limits_{B}\frac{1}{\|\chi_B\|_{L^\varphi}}\int_B |f(x)-f_B|dx<\infty,$$
where $f_B=\frac{1}{|B|}\int_B f(x)dx$ and the supremum is taken over all balls $B$ in $\mathbb R^n$.

Our typical example is $BMO^\varphi(\mathbb R^n)$, called $BMO^{\rm log}(\mathbb R^n)$, related to $\varphi(x,t)=\frac{t}{\log(e+|x|)+ \log(e+t)}$. Clearly, when $\varphi(x,t)\equiv t$, then $BMO^\varphi(\mathbb R^n)$ is just the well-known $BMO(\mathbb R^n)$ of John and Nirenberg. We remark that when $\varphi(x,t)=w(x). t$ with $w\in A_{(n+1)/n}$, then $BMO^\varphi(\mathbb R^n)$ is just $BMO_w(\mathbb R^n)$ was first introduced by Muckenhoupt and Wheeden \cite{MW1, MW2}. There, they proved that $BMO_w(\mathbb R^n)$ is the dual of $H^1_w(\mathbb R^n)$ (see also \cite{Bu}).

\subsection{Quasi-Banach valued sublinear operators}

Recall that a {\sl quasi-Banach space} $\mathcal B$ is a vector space endowed with a quasi-norm $\|\cdot\|_{\mathcal B}$ which is nonnegative, non-degenerate (i.e., $\|f\|_{\mathcal B} = 0$ if and only if $f = 0$), homogeneous, and obeys the quasi-triangle inequality, i.e., there exists a positive constant $\kappa$ no less than $1$ such that for all $f,g\in \mathcal B$, we have $\|f+g\|_{\mathcal B}\leq \kappa(\|f\|_{\mathcal B}+ \|g\|_{\mathcal B})$.

\begin{Definition}\label{quasi-Banach}
Let $\gamma\in (0,1]$. A quasi-Banach space $\mathcal B_\gamma$ with the quasi-norm $\|\cdot\|_{\mathcal B_\gamma}$ is said to be a $\gamma$-quasi-Banach space if  there exists a positive constant $\kappa$ no less than $1$ such that for all $f_j\in \mathcal B_\gamma, j=1,2,..., m$, we have 
$$\Big\|\sum_{j=1}^m f_j\Big\|_{\mathcal B_\gamma}^\gamma\leq \kappa \sum_{j=1}^m \|f_j\|_{\mathcal B_\gamma}^\gamma.$$ 
\end{Definition}

Notice that any Banach space is a $1$-quasi-Banach space, and the quasi-Banach spaces $\ell^p, L^p_w(\mathbb R^n)$ and $H^p_w(\mathbb R^n)$ with $p\in (0,1]$ are typical $p$-quasi-Banach spaces. Also, when $\varphi$ is of uniformly lower type $p\in (0,1]$, the space $H^\varphi(\mathbb R^n)$ is a $p$-quasi-Banach space.

For any given $\gamma$-quasi-Banach space $\mathcal B_\gamma$ with $\gamma\in (0,1]$ and a linear space $\mathcal Y$, an operator $T$ from $\mathcal Y$ to $\mathcal B_\gamma$ is called $\mathcal B_\gamma$-sublinear if there exists a positive constant $\kappa$ no less than $1$ such that 

i) $\|T(f)- T(g)\|_{\mathcal B_\gamma} \leq \kappa \|T(f-g)\|_{\mathcal B_\gamma}$,

ii) for all $f_j\in\mathcal Y, \lambda_j\in \mathbb C,  j=1,...,m$,  we have
$$\Big\|T\Big(\sum_{j=1}^m \lambda_j f_j\Big)\Big\|_{\mathcal B_\gamma}^\gamma \leq \kappa \sum_{j=1}^m |\lambda_j|^\gamma\|T(f_j)\|_{\mathcal B_\gamma}^\gamma.$$

We remark that  if $T$ is linear, then $T$ is $\mathcal B_\gamma$-sublinear. We should point out that if the constant $\kappa$, in Definition \ref{quasi-Banach}, equal 1, then we obtain the notion of $\gamma$-quasi-Banach spaces  introduced in \cite{YZ} (see also \cite{BLYZ}).

\section{Statement of the results}

Our main theorems are the following.

\begin{Theorem}\label{atomic decomposition for new Hardy spaces}
Let $(\varphi, q, s)$ be admissible. Then $H^\varphi(\mathbb R^n)= H^{\varphi, q, s}_{\rm at}(\mathbb R^n)$ with equivalent norms.
\end{Theorem}

Denote by $L^\infty_0(\mathbb R^n)$ the set of all bounded functions  with compact support and zero average. As a consequence of Theorem \ref{atomic decomposition for new Hardy spaces}, we have the following.

\begin{Lemma}
Let $\varphi$ be a growth function satisfying $n q(\varphi)< (n+1) i(\varphi)$. Then, $L^\infty_0(\mathbb R^n)$ is dense in $H^\varphi(\mathbb R^n)$.
\end{Lemma}

We now can present our dual theorem as follows

\begin{Theorem}\label{the dual theorem for new Hardy spaces}
Let $\varphi$ be a growth function satisfying $n q(\varphi)< (n+1) i(\varphi)$. Then, the dual space of $H^\varphi(\mathbb R^n)$ is $BMO^\varphi(\mathbb R^n)$ in the following sense

i)  Suppose $\mathfrak b\in BMO^\varphi(\mathbb R^n)$. Then the linear functional $L_{\mathfrak b}: f\to L_{\mathfrak b}(f):=\int_{\mathbb R^n} f(x)\mathfrak b(x)dx$, initially defined for $L^\infty_0(\mathbb R^n)$, has a bounded extension to $H^\varphi(\mathbb R^n)$.

ii) Conversely, every continuous linear functional on $H^\varphi(\mathbb R^n)$ arises as the above with a unique element $\mathfrak b$ of $BMO^\varphi(\mathbb R^n)$. Moreover $\|\mathfrak b\|_{BMO^\varphi}\approx \|L_{\mathfrak b}\|_{(H^\varphi)^*}$.
\end{Theorem}

Next  result concerns the class of pointwise multipliers for $BMO(\mathbb R^n)$.

\begin{Theorem}\label{a theorem on the class of pointwise multipliers for BMO}
The class of pointwise multipliers for $BMO(\mathbb R^n)$ is  the dual of $L^1(\mathbb R^n)+ H^{\rm log}(\mathbb R^n)$ where $ H^{\rm log}(\mathbb R^n)$ is a Hardy space of Musielak-Orlicz type related to the Musielak-Orlicz function $\theta(x,t)=\frac{t}{\log(e+ |x|)+ \log(e+t)}$.
\end{Theorem}

In order to obtain the finite atomic decomposition, we need the notion of {\sl uniformly locally dominated convergence condition}. A growth function $\varphi$ is said to be satisfy  {\sl uniformly locally dominated convergence condition} if the following holds:

Given $K$ compact set in $\mathbb R^n$. Let $\{f_m\}_{m\geq 1}$ be a sequence of measurable functions  s.t $f_m(x)$ tends to $f(x)$ a.e $x\in\mathbb R^n$. If there exists a nonnegative measurable function $g$ s.t $|f_m(x)|\leq g(x)$ a. e. $x\in \mathbb R^n$ and $\sup_{t>0}\int_K g(x)\frac{\varphi(x,t)}{\int_K \varphi(y,t)dy}dx<\infty$, then $\sup_{t>0}\int_K |f_m(x)-f(x)|\frac{\varphi(x,t)}{\int_K \varphi(y,t)dy}dx$ tends 0.

We remark that the growth functions $\varphi(x,t)=w(x)\Phi(t)$ and $\varphi(x,t)=\frac{t^p}{(\log(e+|x|)+ \log(e+t^p))^p}$, $0<p\leq 1$, satisfy the {\sl uniformly locally dominated convergence condition}.

\begin{Theorem}\label{finite decomposition for new Hardy spaces}
Let  $\varphi$ be a growth function satisfying  uniformly locally dominated convergence condition, and  $(\varphi,q,s)$ be an admissible triplet.

i) If $q\in (q(\varphi), \infty)$ then $\|\cdot\|_{H^{\varphi,q,s}_{\rm fin}}$ and $\|\cdot\|_{H^\varphi}$ are equivalent quasi-norms on  $H^{\varphi,q,s}_{\rm fin}(\mathbb R^n)$.

ii) $\|\cdot\|_{H^{\varphi,\infty,s}_{\rm fin}}$ and $\|\cdot\|_{H^\varphi}$ are equivalent quasi-norms on  $H^{\varphi,\infty,s}_{\rm fin}(\mathbb R^n)\cap C(\mathbb R^n)$.
\end{Theorem}

As an application, we obtain criterions for boundedness of quasi-Banach valued sublinear operators in $H^\varphi(\mathbb R^n)$.

\begin{Theorem}\label{the boundedness of sublinear operators on new Hardy spaces}
Let  $\varphi$ be a growth function satisfying  uniformly locally dominated convergence condition, $(\varphi,q,s)$ be an admissible triplet, $\varphi$ be of uniformly upper type $\gamma\in (0,1]$, and $\mathcal B_\gamma$ be a quasi-Banach space. Suppose  one of the following holds:

i) $q\in (q(\varphi),\infty)$, and $T: H^{\varphi,q,s}_{\rm fin}(\mathbb R^n)\to \mathcal B_\gamma$ is a $\mathcal B_\gamma$-sublinear operator such that
$$A= \sup\{\|Ta\|_{\mathcal B_\gamma}: a\;\mbox{is a}\; (\varphi,q,s){\rm -atom}\}<\infty;$$

ii) T is a $\mathcal B_\gamma$-sublinear operator defined on continuous $(\varphi,\infty,s)$-atoms such that
$$A= \sup\{\|Ta\|_{\mathcal B_\gamma}: a\;\mbox{is a continuous}\; (\varphi,\infty,s){\rm -atom}\}<\infty.$$
Then there exists a unique bounded $\mathcal B_\gamma$-sublinear operator $\widetilde T$ from $H^\varphi(\mathbb R^n)$ to $\mathcal B_\gamma$ which extends $T$.
\end{Theorem}

\section{Some basic lemmas on growth functions}

We start by the following lemma.

\begin{Lemma}\label{the basis lemma 1}
i) Let $\varphi$ be a growth function. Then $\varphi$ is uniformly $\sigma$-quasi-subadditive on $\mathbb R^n\times [0,\infty)$, i.e. there exists a constant $C>0$ such that
$$\varphi(x,\sum_{j=1}^\infty t_j)\leq C \sum_{j=1}^\infty \varphi(x,t_j),$$
for all $(x,t_j)\in\mathbb R^n\times [0,\infty)$, $j=1,2,...$

ii) Let $\varphi$ be a growth function and $\widetilde\varphi(x,t):= \int_0^t \frac{\varphi(x,s)}{s}ds$ for $(x,t)\in\mathbb R^n\times [0,\infty)$. Then $\widetilde\varphi$ is a growth function equivalent to $\varphi$, moreover, $\widetilde\varphi(x,\cdot)$ is continuous and strictly increasing.

iii) A Musielak-Orlicz function $\varphi$ is a  growth function if and only if $\varphi$ is of positive uniformly lower type and uniformly quasi-concave, i.e. there exists a constant $C>0$ such that
$$\lambda\varphi(x,t)+ (1-\lambda)\varphi(x,s)\leq C \varphi(x, \lambda t+ (1-\lambda)s),$$
for all $x\in\mathbb R^n, t,s\in [0,\infty)$ and $\lambda\in [0,1]$.
\end{Lemma}
\begin{proof}
 i) We just need to consider the case when $\sum_{j=1}^\infty t_j>0$. Then it follows from the fact that
$$\frac{t_k}{\sum_{j=1}^\infty t_j}\varphi(x,\sum_{j=1}^\infty t_j)\leq C \varphi(x, t_k)$$
by $\varphi$ is of uniformly upper type 1.

ii) Since $\varphi$ is a growth function, it is easy to see that $\widetilde\varphi(x,\cdot)$ is continuous and strictly increasing. Moreover, there exists $p>0$ such that $\varphi$ is of uniformly lower type $p$. Hence,
\begin{equation}\label{the basis lemma 1, 1}
\widetilde\varphi(x,t)= \int_0^t \frac{\varphi(x,s)}{s}ds\leq C\frac{\varphi(x,t)}{t^p}\int_0^t \frac{1}{s^{1-p}}ds\leq C \varphi(x,t).
\end{equation}

On the other hand, since $\varphi$ is of uniformly upper type 1, we get
\begin{equation}\label{the basis lemma 1, 2}
\widetilde\varphi(x,t)= \int_0^t \frac{\varphi(x,s)}{s}ds\geq C^{-1}\int_0^t \frac{\varphi(x,t)}{t}ds\geq C^{-1}\varphi(x,t).
\end{equation}

Combining (\ref{the basis lemma 1, 1}) and (\ref{the basis lemma 1, 2}), we obtain $\widetilde\varphi\approx \varphi$, and thus $\widetilde\varphi$ is a growth function.

iii) Suppose  $\varphi$ is a growth function. By (ii), $\varphi$ is equivalent to $\widetilde{\widetilde\varphi}$. On the other hand, $\frac{\partial\widetilde{\widetilde\varphi}}{\partial t}(x,t)=\frac{\widetilde\varphi(x,t)}{t}$ is uniformly quasi-decreasing in $t$. Hence, $\widetilde{\widetilde\varphi}$ is uniformly quasi-concave, and thus is $\varphi$.

The converse is easy by taking $s=0$. We omit the details.
\end{proof}

\begin{Remark}
Let us observe that the results stated in Section 3 are invariant under change of equivalent growth functions. By Lemma \ref{the basis lemma 1}, in the future, we always consider a  growth function $\varphi$ of positive uniformly lower type, of  uniformly upper type 1 (or, equivalently,  uniformly quasi-concave), and so that $\varphi(x,\cdot)$ is continuous and strictly increasing for all $x\in\mathbb R^n$.
\end{Remark}

\begin{Lemma}\label{the basis lemma 2}
Let $\varphi$ be a growth function. Then

i) $\displaystyle\int_{\mathbb R^n}\varphi\Big(x, \frac{|f(x)|}{\|f\|_{L^\varphi}}\Big)dx=1$ for all $f\in L^\varphi(\mathbb R^n)\setminus \{0\}$.

ii) $\lim_{k\to\infty}\|f_k\|_{L^\varphi}=0$ if and only if $\lim_{k\to\infty}\int_{\mathbb R^n}\varphi(x,|f_k(x)|)dx=0$.
\end{Lemma}
\begin{proof}
 Statement  $(i)$ follows from the fact that the function
$$\vartheta(t):= \int_{\mathbb R^n}\varphi(x, t|f(x)|)dx,$$
$t\in [0,\infty)$, is continuous by the dominated convergence theorem since $\varphi(x,\cdot)$ is continuous.

Statement $(ii)$ follows from the fact that
$$\|f\|_{L^\varphi}\leq C \max\Big\{\int_{\mathbb R^n}\varphi(x,|f(x)|)dx, \Big(\int_{\mathbb R^n}\varphi(x,|f(x)|)dx\Big)^{1/p}\Big\},$$
and 
$$\int_{\mathbb R^n}\varphi(x,|f(x)|)dx\leq C\max\Big\{\|f\|_{L^\varphi}, (\|f\|_{L^\varphi})^p\Big\}$$
for some $p\in (0,i(\varphi))$. 
\end{proof}

\begin{Lemma}\label{the basis lemma 3}
Given $c$ is a positive constant. Then, there exists a constant $C>0$ such that

i) The inequality $\int_{\mathbb R^n}\varphi\Big(x,\frac{|f(x)|}{\lambda}\Big)dx\leq c$,
for $\lambda>0$, implies
$$\|f\|_{L^\varphi}\leq C \lambda.$$

ii) The inequality $\sum_j \varphi\Big(B_j, \frac{t_j}{\lambda}\Big)\leq c$, for $\lambda>0$, implies
$$\inf\Big\{\alpha>0: \sum_j \varphi\Big(B_j, \frac{t_j}{\alpha}\Big)\leq 1\Big\}\leq C \lambda.$$
\end{Lemma}

\begin{proof}
The proofs are simple since we may take $C=(1+c.C_p)^{1/p}$, for some $p\in (0,i(\varphi))$, where  $C_p$ is such that (\ref{lower type}) holds.
\end{proof}

\begin{Lemma}\label{the basis lemma 4}
Let $(\varphi,q,s)$ be an admissible triplet. Then there exists a positive constant $C$ such that 
$$\sum_{j=1}^\infty \|b_j\|_{L^q_\varphi(B_j)}\|\chi_{B_j}\|_{L^\varphi}\leq C \Lambda_q(\{b_j\}),$$
for all $f=\sum_{j=1}^\infty b_j\in H^{\varphi,q,s}_{at}(\mathbb R^n)$ where $b_j$'s are multiples of $(\varphi,q,s)$-atoms supported in balls $B_j$'s.
\end{Lemma}
\begin{proof}
Since $\varphi$ is of uniformly upper type 1,  there exists a positive constant $c>0$ such that
$$\varphi\Big(x, \frac{\|b_i\|_{L^q_\varphi(B_i)}}{\sum_{j=1}^\infty \|b_j\|_{L^q_\varphi(B_j)}\|\chi_{B_j}\|_{L^\varphi}}\Big)\geq c \frac{\|b_i\|_{L^q_\varphi(B_i)}\|\chi_{B_i}\|_{L^\varphi}}{\sum_{j=1}^\infty \|b_j\|_{L^q_\varphi(B_j)}\|\chi_{B_j}\|_{L^\varphi}}\varphi\Big(x, \frac{1}{\|\chi_{B_i}\|_{L^\varphi}}\Big)$$
for all $x\in\mathbb R^n, i\geq 1$. Hence, for all $i\geq 1$,
$$\varphi\Big(B_i, \frac{\|b_i\|_{L^q_\varphi(B_i)}}{\sum_{j=1}^\infty \|b_j\|_{L^q_\varphi(B_j)}\|\chi_{B_j}\|_{L^\varphi}}\Big)\geq c \frac{\|b_i\|_{L^q_\varphi(B_i)}\|\chi_{B_i}\|_{L^\varphi}}{\sum_{j=1}^\infty \|b_j\|_{L^q_\varphi(B_j)}\|\chi_{B_j}\|_{L^\varphi}}$$
since $\int_{B_i}\varphi\Big(x, \frac{1}{\|\chi_{B_i}\|_{L^\varphi}}\Big)dx=1$ by Lemma \ref{the basis lemma 2}. It follows that
$$\sum_{i=1}^\infty \varphi\Big(B_i, \frac{\|b_i\|_{L^q_\varphi(B_i)}}{\sum_{j=1}^\infty \|b_j\|_{L^q_\varphi(B_j)}\|\chi_{B_j}\|_{L^\varphi}}\Big)\geq c.$$
We  deduce from the above that
$$\sum_{j=1}^\infty \|b_j\|_{L^q_\varphi(B_j)}\|\chi_{B_j}\|_{L^\varphi}\leq C \Lambda_q(\{b_j\}),$$
with $C=(1+ C_p/c )^{1/p}$ for some $p\in (0,i(\varphi))$, where $C_p$ is such that (\ref{lower type}) holds.
\end{proof}

\begin{Lemma}\label{the basis lemma 5}
Let $\varphi\in \mathbb A_q, 1<q<\infty$. Then, there exists a positive constant $C$ such that

i) For all ball $B(x_0,r), \lambda>1$, and $t\in [0,\infty)$, we have
$$\varphi(B(x_0, \lambda r), t)\leq C \lambda^{nq} \varphi(B(x_0,r), t).$$

ii) For all ball $B(x_0,r)$ and $t\in [0,\infty)$, we have
$$\int_{B^c}\frac{\varphi(x,t)}{|x-x_0|^{nq}}dx\leq C \frac{\varphi(B,t)}{r^{nq}}.$$

iii) For all ball $B$, $f$ measurable and $t\in (0,\infty)$, we have
$$\Big(\frac{1}{|B|}\int_B |f(x)| dx\Big)^q\leq C \frac{1}{\varphi(B,t)}\int_B |f(x)|^q \varphi(x,t)dx.$$

iv) For all $f$ measurable and $t\in [0,\infty)$, we have
$$\int_{\mathbb R^n}\mathcal Mf(x)^q \varphi(x,t)dx\leq C\int_{\mathbb R^n}|f(x)|^q \varphi(x,t)dx,$$
where $\mathcal M$ is the classical  Hardy-Littlewood maximal operator defined by
$$\mathcal Mf(x)= \sup\limits_{x\in B\rm -ball}\frac{1}{|B|}\int_B |f(y)|dy\;,\quad x\in\mathbb R^n.$$
\end{Lemma}

 In the setting  $\varphi(x,t)=w(x)\Phi(t)$, $w\in A_\infty$ and $\Phi$ a Orlicz function, the above lemma is well-known as a classical result in the theory of Muckenhoupt weight (see \cite{GR}).  Since $\varphi$ satisfies uniformly Muckenhoupt condition, the proof of Lemma \ref{the basis lemma 5} is a slight modification of the classical result. We omit the details.

\section{Atomic decompositions}

The  purpose of this section is prove the atomic decomposition theorem (Theorem \ref{atomic decomposition for new Hardy spaces}). The construction is by now standard, but the estimates require the preliminary lemmas. For the reader convenience, we give all steps of the proof, even if only the generalization to our framework is new. 

We first  introduce a class of Hardy spaces containing  the Hardy space of Musielak-Orlicz type $H^\varphi(\mathbb R^n)$ as a particular case.

\begin{Definition}
For $m\in\mathbb N$, we denote by $H^\varphi_m(\mathbb R^n)$ the space of all distributions $f$ such that $f^*_m\in L^\varphi(\mathbb R^n)$ with the (quasi-)norm
$$\|f\|_{H^\varphi_m}:= \|f^*_m\|_{L^\varphi}.$$
\end{Definition}
Clearly, $H^\varphi(\mathbb R^n)$ is a special case associated with $m=m(\varphi)$.

\subsection{Some basic properties concerning $H^\varphi_m(\mathbb R^n)$ and $H^{\varphi,q,s}_{at}(\mathbb R^n)$}
We start by the following proposition.

\begin{Proposition}\label{the inclusion of Hardy spaces into the Schwartz space}
For $m\in\mathbb N$, we have $ H^\varphi_m(\mathbb R^n)\subset \mathcal S'(\mathbb R^n)$ and the inclusion is continuous.
\end{Proposition}
\begin{proof}
 Let $f\in  H^\varphi_m(\mathbb R^n)$. For any $\phi\in \mathcal S(\mathbb R^n)$, and $x\in B(0,1)$, we write 
$$\left\langle{f,\phi}\right\rangle=f*\widetilde \phi (0)=f*\psi(x),$$
where $\psi(y)=\widetilde \phi(y-x)=\phi(x-y)$ for all $y\in\mathbb R^n$.
 
It is easy to verify that $\sup\limits_{x\in B(0,1), y\in \mathbb R^n}\frac{1+|y|}{1+|y-x|}\leq 2$. Consequently,
\begin{eqnarray*}
|\left\langle{f,\phi}\right\rangle|=|f*\psi(x)|
&\leq&
2^{(m+2)(n+1)}\|\phi\|_{\mathcal S_m}\inf\limits_{x\in B(0,1)} f^*_m(x)\\
&\leq&
2^{(m+2)(n+1)}\|\phi\|_{\mathcal S_m}\|\chi_{B(0,1)}\|^{-1}_{L^\varphi}\|f\|_{ H^\varphi_m}.
\end{eqnarray*}
This implies that $f\in\mathcal S'(\mathbb R^n)$ and the inclusion is continuous. 
\end{proof}

The following proposition gives the completeness of $ H^\varphi_m(\mathbb R^n)$.

\begin{Proposition}\label{the completeness of new Hardy spaces}
The space $ H^\varphi_m(\mathbb R^n)$ is complete.
\end{Proposition}
\begin{proof}
 In order to prove the completeness of $ H^\varphi_m(\mathbb R^n)$, it suffices to prove that for every sequence $\{f_j\}_{j\geq 1}$ with $\|f_j\|_{ H^\varphi_m}\leq 2^{-j}$ for any $j\geq 1$, the series $\sum_j f_j$ converges in $ H^\varphi_m(\mathbb R^n)$. Let us now take $p>0$ such that $\varphi$ is of uniformly lower type $p$. Then, for any $j\geq 1$, 
\begin{equation}\label{the completeness of new Hardy spaces 1}
\int_{\mathbb R^n}\varphi(x, (f_j)^*_m(x))dx\leq C (2^{-j})^p\int_{\mathbb R^n} \varphi\Big(x, \frac{(f_j)^*_m(x)}{2^{-j}}\Big)dx\leq C 2^{-jp}.
\end{equation}

Since $\{\sum_{i=1}^j f_i\}_{j\geq 1}$ is a Cauchy sequence in $ H^\varphi_m(\mathbb R^n)$, by  Proposition \ref{the inclusion of Hardy spaces into the Schwartz space} and the completeness of $\mathcal S'(\mathbb R^n)$, $\{\sum_{i=1}^j f_i\}_{j\geq 1}$ is also a Cauchy sequence in $\mathcal S'(\mathbb R^n)$ and thus converges to some $f\in\mathcal S'(\mathbb R^n)$. This implies that, for every $\phi\in \mathcal S(\mathbb R^n)$, the series $\sum_j f_j*\phi$ converges to $f*\phi$ pointwisely. Therefore $f^*_m(x)\leq \sum_j (f_j)^*_m(x)$ and $(f-\sum_{j=1}^k f_j)^*_m(x)\leq \sum_{j\geq {k+1}} (f_j)^*_m(x)$ for all $x\in\mathbb R^n, k\geq 1$. Combining this and (\ref{the completeness of new Hardy spaces 1}), we obtain
\begin{eqnarray*}
\int_{\mathbb R^n}\varphi(x,(f-\sum_{j=1}^kf_j)^*_m(x))dx
&\leq&
C \sum_{j\geq {k+1}} \int_{\mathbb R^n}\varphi(x, (f_j)^*_m(x))dx\\
&\leq&
C \sum_{j\geq {k+1}} 2^{-jp}\to 0, 
\end{eqnarray*}
as $k\to\infty$, here we used  Lemma \ref{the basis lemma 1}. Thus, the series $\sum_j f_j$ converges to $f$ in $ H^\varphi_m(\mathbb R^n)$ by Lemma \ref{the basis lemma 2}.
This completes the proof.
\end{proof}

\begin{Corollary}
The Hardy space of Musielak-Orlicz type $ H^\varphi(\mathbb R^n)$ is complete.
\end{Corollary}

The following lemma and its corollary show that $(\varphi, q,s)$-atoms are in $H^\varphi(\mathbb R^n)$. Furthermore, it is the necessary estimate for  proving that  $H^{\varphi,q,s}_{\rm at}(\mathbb R^n)\subset H^\varphi(\mathbb R^n)$ and the inclusion is continuous, see Theorem \ref{part 1 for the atomic decomposition} below.

\begin{Lemma}\label{a basis estimate for atoms in new Hardy spaces}
Let $(\varphi,q,s)$ be an admissible triplet and $m\geq s$. Then, there exists a constant $C=C(\varphi, q, s, m)$ such that
$$\int_{\mathbb R^n}\varphi(x, f^*_m(x))dx\leq C \varphi(B, \|f\|_{L^q_\varphi(B)}),$$
for all $f$ multiples of $(\varphi,q,s)$-atom associated with ball $B=B(x_0,r)$. 
\end{Lemma}
\begin{proof}
 The case $q=\infty$ is easy and will be omitted. We just consider $q\in(q(\varphi),\infty)$. Now let us set $\widetilde B=B(x_0,9r)$, and write
\begin{eqnarray*}
\int_{\mathbb R^n}\varphi(x, f^*_m(x))dx
&=&
\int_{\widetilde B}\varphi(x,f^*_m(x))dx+ \int_{(\widetilde B)^c}\varphi(x,f^*_m(x))dx\\
&=&
I+II.
\end{eqnarray*}
Since $\varphi$ is of uniformly upper type 1, by H\"older inequality, we get
\begin{eqnarray*}
I
&=&
\int_{\widetilde B}\varphi(x, f^*_m(x))dx\leq C \int_{\widetilde B}\Big(\frac{f^*_m(x)}{\|f\|_{L^q_\varphi(B)}}+ 1\Big)\varphi(x,\|f\|_{L^q_\varphi(B)})dx\\
&\leq&
C \varphi(\widetilde B,\|f\|_{L^q_\varphi(B)})\\
&+&
C \frac{1}{\|f\|_{L^q_\varphi(B)}}\Big(\int_{\widetilde B}|f^*_m(x)|^q \varphi(x,\|f\|_{L^q_\varphi(B)})dx\Big)^{1/q} \varphi(\widetilde B, \|f\|_{L^q_\varphi(B)})^{(q-1)/q}\\
&\leq&
C \varphi(B,\|f\|_{L^q_\varphi(B)})+ C \frac{1}{\|f\|_{L^q_\varphi(B)}}\|f\|_{L^q_\varphi(\widetilde B)} \varphi(\widetilde B, \|f\|_{L^q_\varphi(B)})\\
&\leq&
C \varphi(B,\|f\|_{L^q_\varphi(B)}).
\end{eqnarray*}
We used the fact $f^*_m(x)\leq C(m)\mathcal Mf(x)$  and Lemma \ref{the basis lemma 5}.

To estimate $II$, we note that since $m\geq s$, there exists a constant $C=C(m)$ such that
$$\Big|\phi\Big(\frac{x-y}{t}\Big)-\sum_{|\alpha|\leq s}\frac{\partial^\alpha\phi(\frac{x-x_0}{t})}{\alpha !}\Big(\frac{x_0-y}{t}\Big)^\alpha\Big|\leq C t^n\frac{|y-x_0|^{s+1}}{|x-x_0|^{n+s+1}}$$
for all $\phi\in\mathcal S_m(\mathbb R^n), t>0, x\in (\widetilde B)^c, y\in B$. Therefore
\begin{eqnarray*}
|f*\phi_t(x)|
&=&
\frac{1}{t^n}\Big|\int_{B}f(y)\Big[\phi\Big(\frac{x-y}{t}\Big)-\sum_{|\alpha|\leq s}\frac{\partial^\alpha\phi(\frac{x-x_0}{t})}{\alpha !}\Big(\frac{x_0-y}{t}\Big)^\alpha\Big]dy\Big|\\
&\leq&
C\int_{B}|f(y)|\frac{|y-x_0|^{s+1}}{|x-x_0|^{n+s+1}}dy\\
&\leq&
C\frac{r^{s+1}}{|x-x_0|^{n+s+1}}\Big(\int_{B}|f(y)|^q \varphi(y,\lambda)dy\Big)^{1/q}\Big(\int_{B}[\varphi(y,\lambda)]^{-1/(q-1)}dy\Big)^{(q-1)/q}\\
&\leq&
C \|f\|_{L^q_\varphi(B)}\Big(\frac{r}{|x-x_0|}\Big)^{n+s+1}.
\end{eqnarray*}
For any $\lambda>0$,  we used that $\int_B \varphi(y,\lambda)dy(\int_B [\varphi(y,\lambda)]^{-1/(q-1)}dy)^{q-1}\leq C |B|^q$ since $\varphi\in \mathbb A_q$. As a consequence, we get
$$f^*_m(x)\leq C(m) \sup_{\phi\in\mathcal S_m(\mathbb R^n)}\sup\limits_{t>0}|f*\phi_t(x)|\leq C \|f\|_{L^q_\varphi(B)}\Big(\frac{r}{|x-x_0|}\Big)^{n+s+1}.$$
By $s\geq m(\varphi)$, there exists $p\in (0,i(\varphi))$ such that $(n+s+1)p> nq(\varphi)$. Hence, by Lemma \ref{the basis lemma 5},
\begin{eqnarray*}
II=\int_{(\widetilde B)^c}\varphi(x, f^*_m(x))dx
&\leq&
C \int_{(\widetilde B)^c}\Big(\frac{r}{|x-x_0|}\Big)^{(n+s+1)p}\varphi(x,\|f\|_{L^q_\varphi(B)})dx\\
&\leq&
C r^{(n+s+1)p}\frac{\varphi(\widetilde B, \|f\|_{L^q_\varphi(B)})}{(9r)^{(n+s+1)p}}\\
&\leq&
C \varphi(B, \|f\|_{L^q_\varphi(B)}).
\end{eqnarray*}
This completes the proof.
\end{proof}

\begin{Corollary}
There exists a constant $C=C(\varphi, q, s)>0$ such that
$$\|a\|_{H^\varphi}\leq C,$$
for all $(\varphi, q, s)$-atom $a$.
\end{Corollary}

\begin{Theorem}\label{part 1 for the atomic decomposition}
Let $(\varphi,q,s)$ be an admissible triplet and $m\geq s$. Then
$$ H^{\varphi,q,s}_{\rm at}(\mathbb R^n)\subset H^\varphi_m(\mathbb R^n),$$
moreover, the inclusion is continuous.
\end{Theorem}
\begin{proof}
 For any $0\ne f\in  H^{\varphi,q,s}_{\rm at}(\mathbb R^n)$. Let $f=\sum_j b_j$ be an atomic decomposition of $f$, with supp $b_j\subset B_j$, $j=1,2,...$ For all $\phi\in\mathcal S(\mathbb R^n)$, the series $\sum_j b_j*\phi$ converges to $f*\phi$ pointwise since $f=\sum_j b_j$ in $\mathcal S'$. Hence $f^*_m(x)\leq \sum_{j} (b_j)^*_m(x)$. By applying Lemma \ref{a basis estimate for atoms in new Hardy spaces}, we obtain
\begin{eqnarray*}
\int_{\mathbb R^n}\varphi\Big(x, \frac{f^*_m(x)}{\Lambda_q(\{b_j\})}\Big)dx
&\leq&
C \sum_{j} \int_{\mathbb R^n}\varphi\Big(x, \frac{(b_j)^*_m(x)}{\Lambda_q(\{b_j\})}\Big)dx\\
&\leq&
C \sum_{j} \varphi\Big(B_j, \frac{\|b_j\|_{L^q_\varphi(B_j)}}{\Lambda_q(\{b_j\})}\Big)\\
&\leq&
C.
\end{eqnarray*}

This implies that $\|f\|_{H^\varphi_m}\leq C \Lambda_q(\{b_j\})$ (see Lemma \ref{the basis lemma 3}) for any atomic decomposition $f=\sum_j b_j$, and thus, $\|f\|_{H^\varphi_m}\leq C \|f\|_{ H^{\varphi,q,s}_{\rm at}}$. 
\end{proof}

\subsection{Calder\'on-Zygmund decompositions}\label{subsection for Calderon-Zygmund decompositions}

Throughout this subsection, we fix $m$ and $s$ so that $ m, s \geq m(\varphi)$. For a given $\lambda>0$, we set $\Omega=\{x\in\mathbb R^n: f^*_m(x)>\lambda\}$. Observe that $\Omega$ is open. Hence by Whitney's lemma, there exist $x_1,x_2,...$ in $\Omega$ and $r_1,r_2,...>0$ such that

(i) $\Omega=\cup_j B(x_j,r_j)$,

(ii) the balls $B(x_j,r_j/4)$, $j=1,2,...$, are disjoint,

(iii) $B(x_j,18 r_j)\cap \Omega^c=\emptyset$, but $B(x_j,54 r_j)\cap \Omega^c\ne \emptyset$, for any $j=1,2,...,$

(iv) there exists $L\in \mathbb N$ (depending only on $n$) such that no point of $\Omega$ lies in more than $L$ of the balls $B(x_j, 18r_j)$, $j=1,2,...$

We fix once for all, a function $\theta\in C_0^\infty(\mathbb R^n)$ such that supp $\theta\subset B(0,2)$, $0\leq \theta\leq 1$, $\theta=1$ on $B(0,1)$, and set $\theta_j(x)=\theta((x-x_j)/r_j)$, for j=1,2,...  Obviously, supp $\theta_j\subset B(x_j,2r_j)$, $j=1,2,...$, and $1\leq \sum_j \theta_j\leq L$ for all $x\in \Omega$. Hence if we set $\zeta_j(x)=\theta_j(x)/\sum_{i=1}^\infty \theta_i(x)$ if $x\in\Omega$ and $\zeta_j(x)=0$ if $x\in \Omega^c$, $j=1,2,...$, then supp $\zeta_j\subset B(x_j,2r_j)$, $0\leq \zeta_j\leq 1$, $\sum_j\zeta_j=\chi_\Omega$, and $L^{-1}\leq \zeta_j\leq 1$ on $B(x_j,r_j)$. The family $\{\zeta_j\}_j$ forms a smooth partition of unity of $\Omega$. let $s\in \mathbb N$ be some fixed natural number and $\mathcal P_s(\mathbb R^n)$ (or simply $\mathcal P_s$) denote the linear space of polynomials in $n$ variables of degree less than $s$. For each $j$, we consider the inner product $\left\langle{P,Q}\right\rangle_{j}=\frac{1}{\int_{\mathbb R^n}\zeta_j(x)dx}\int_{\mathbb R^n}P(x)Q(x) \zeta_j(x)dx$ for $P,Q\in \mathcal P_s$. Then $(\mathcal P_s,\left\langle{\cdot,\cdot}\right\rangle_j)$ is a finite dimensional Hilbert space. Let $f\in \mathcal S'$. Since $f$ induces a linear functional on $\mathcal P_s$ via $Q\to \frac{1}{\int_{\mathbb R^n}\zeta_j(x)dx}\int_{\mathbb R^n}f(x)Q(x) \zeta_j(x)dx$, by the Riesz theorem, there exists a unique polynomial $P_j\in \mathcal P_s$ such that for all $Q\in \mathcal P_s$, $\left\langle{P_j,Q}\right\rangle_j=\frac{1}{\int_{\mathbb R^n}\zeta_j(x)dx}\int_{\mathbb R^n}f(x)Q(x) \zeta_j(x)dx.$ For each $j$, $j=1,2,...$, we define $b_j=(f-P_j)\zeta_j$, and note $B_j=B(x_j,r_j)$, $\widetilde B_j=B(x_j,9r_j)$. Then, it is easy to see that $\int_{\mathbb R^n}b_j(x)Q(x)dx=0$ for all $Q\in\mathcal P_s$. It turns out, in the case of interest, that the series $\sum_j b_j$ converges in $\mathcal S'$. In this case, we set $g=f-\sum_j b_j$, and we call the representation $f=g+\sum_j b_j$ a Calder\'on-Zygmund decomposition of $f$ of degree $s$ and height $\lambda$ associated to $f^*_m$. 

For any $j=1,2,...$, we denote $B_j=B(x_j,r_j)$ and $\widetilde B_j=B(x_j,9r_j)$. Then we have the following lemma which proof can be found in \cite[Chapter 3]{FoS}.

\begin{lemmaa}
There are four constant $c_1,c_2,c_3,c_4$, independent of $f,j,$ and $\lambda$, such that

i)
$$\sup\limits_{|\alpha|\leq N,x\in \mathbb R^n} r_j^{|\alpha|}|\partial^\alpha \zeta_j(x)|\leq c_1.$$

ii) 
$$\sup\limits_{x\in \mathbb R^n}|P_j(x)\zeta_j(x)|\leq c_2 \lambda.$$

iii)
$$(b_j)^*_m(x)\leq c_3 f^*_m(x),\quad\mbox{for all}\;\; x\in \widetilde B_j.$$

iv)
$$(b_j)^*_m(x)\leq c_4 \lambda (r_j/{|x-x_j|})^{n+m_s},\quad\mbox{ for all}\;\; x\notin \widetilde B_j,$$
where $m_s=\min\{s+1,m+1\}$.
\end{lemmaa}

\begin{Lemma}\label{basis lemma for atomic decomposition 1}
For all $f\in  H^\varphi_m(\mathbb R^n)$, there exists a geometric constant $C$, independent of $f,j$, and $\lambda$, such that, 
$$\int_{\mathbb R^n}\varphi\Big(x, (b_j)^*_m(x)\Big)dx\leq C \int_{\widetilde B_j}\varphi(x, f^*_m(x))dx.$$

Moreover, the series $\sum_j b_j$ converges in $ H^\varphi_m(\mathbb R^n)$, and

$$\int_{\mathbb R^n}\varphi\Big(x, (\sum_j b_j)^*_m(x)\Big)dx\leq C\int_{\Omega}\varphi(x, f^*_m(x))dx.$$
\end{Lemma}

\begin{proof}
 As $m,s\geq m(\varphi)$, $m_s=\min\{s+1, m+1\}> n(q(\varphi)/i(\varphi)-1)$. Hence, there exist $q>q(\varphi)$ and $0<p<i(\varphi)$ such that $m_s>n(q/p-1)$, deduce that $(n+m_s)p>nq$. Therefore, $\varphi\in \mathbb A_{(n+m_s)p/n}$ and $\varphi$ is of uniformly lower type $p$. Thus, there exists a positive constant $C$, independent of $f,j,$ and $\lambda$, such that
\begin{eqnarray*}
\int_{(\widetilde B_j)^c}\varphi(x, \lambda (r_j/{|x-x_j|})^{n+m_s})dx
&\leq&
C \int_{(\widetilde B_j)^c}\Big(\frac{r_j}{|x-x_j|}\Big)^{(n+m_s)p}\varphi(x, \lambda)dx\\
&\leq&
C (r_j)^{(n+m_s)p} \frac{\varphi(\widetilde B_j, \lambda)}{(9r_j)^{(n+m_s)p}}\\
&\leq&
C \int_{\widetilde B_j}\varphi(x, f^*_m(x))dx,
\end{eqnarray*}
since $r_j/{|x-x_j|}<1$ and $ f^*_m>\lambda$ on $\widetilde B_j$. Combining this and Lemma A, we get

\begin{eqnarray*}
\int_{\mathbb R^n}\varphi\Big(x, (b_j)^*_m(x)\Big)dx
&\leq&
C\Big[\int_{\widetilde B_j} \varphi(x, f^*_m(x))dx+ \int_{(\widetilde B_j)^c} \varphi(x, \lambda (r_j/|x-x_j|)^{n+m_s})dx\Big]\\
&\leq&
C \int_{\widetilde B_j} \varphi(x, f^*_m(x))dx.
\end{eqnarray*}
As a consequence of the above estimate, since $\sum_j \chi_{\widetilde B_j}\leq L$ and $\Omega=\cup_j \widetilde B_j$, we obtain
\begin{eqnarray*}
\sum_j \int_{\mathbb R^n}\varphi\Big(x, (b_j)^*_m(x)\Big)dx
&\leq&
C\sum_j \int_{\widetilde B_j}\varphi(x, f^*_m(x))dx\\
&\leq&
C\int_{\Omega}\varphi(x, f^*_m(x))dx.
\end{eqnarray*}
This implies that the series $\sum_j b_j$ converges in $ H^\varphi_m(\mathbb R^n)$ by completeness of $ H^\varphi_m(\mathbb R^n)$. Moreover,
$$\int_{\mathbb R^n}\varphi\Big(x,(\sum_j b_j)^*_m(x)\Big)dx\leq C \int_{\Omega}\varphi(x, f^*_m(x))dx.$$
\end{proof}

Let $q\in [1,\infty]$. We denote by $L^q_{\varphi(\cdot,1)}(\mathbb R^n)$ the usually weighted Lebesgue space with the Muckenhoupt weight $\varphi(x,1)$. Then, we have the following.

\begin{lemmab}[see \cite{BLYZ}, Lemma 4.8]
Let $q\in (q(\varphi),\infty]$. Assume that $f\in L^q_{\varphi(\cdot,1)}(\mathbb R^n)$, then the series $\sum_j b_j$ converges in $L^q_{\varphi(\cdot,1)}(\mathbb R^n)$ and there exists a constant $C$, independent of $f,j,$ and $\lambda$ such that $\|\sum_j |b_j|\|_{L^q_{\varphi(\cdot,1)}}\leq C \|f\|_{L^q_{\varphi(\cdot,1)}}$.
\end{lemmab}

\begin{Remark}
By Lemma B, the series $\sum_j |b_j|$, and thus the series $\sum_j b_j$, converges almost everywhere on $\mathbb R^n$.  
\end{Remark}

\begin{lemmac}[see \cite{FoS}, Lemma 3.19]
 Suppose that the series $\sum_j b_j$ converges in $\mathcal S'(\mathbb R^n)$. Then, there exists a positive constant $C$, independent of $f,j,$ and $\lambda$, such that for all $x\in\mathbb R^n$,
$$g^*_m(x)\leq C \lambda \sum_j\Big(\frac{r_j}{|x-x_j|+r_j}\Big)^{n+m_s}+ f^*_m(x)\chi_{\Omega^c}(x).$$
\end{lemmac}

\begin{Lemma}\label{basis lemma for atomic decomposition 2}
For any $q\in (q(\varphi),\infty)$ and $f\in H^\varphi_m(\mathbb R^n)$. Then $g^*_m\in L^q_{\varphi(\cdot,1)}(\mathbb R^n)$, and there exists a positive constant $C$, independent of $f,j,$ and $\lambda$, such that
$$\int_{\mathbb R^n}[g^*_m(x)]^q \varphi(x,1)dx\leq C\lambda^q \max\{1/\lambda, 1/\lambda^p\}\int_{\mathbb R^n}\varphi(x, f^*_m(x))dx.$$
\end{Lemma}
\begin{proof}
 For any $j=1,2,...$ and $x\in\mathbb R^n$, we have
$$\Big(\frac{r_j}{|x-x_j|+r_j}\Big)^n=\frac{1}{|B(x_j,|x-x_j|+r_j)|}\int_{B(x_j,|x-x_j|+r_j)}\chi_{B_j}(y)dy\leq \mathcal M(\chi_{B_j})(x)$$
since $B_j\subset B(x_j,|x-x_j|+r_j)$. 

Therefore, by $L^{rq}_{\varphi(\cdot,1)}$-boundedness of vector-valued maximal functions (see \cite[Theorem 3.1]{AJ}), where $r:=(n+ m_s)/n>1$, we obtain that
\begin{eqnarray*}
\int_{\mathbb R^n}\Big[\sum_j\Big(\frac{r_j}{|x-x_j|+r_j}\big)^{n+m_s}\Big]^q \varphi(x,1)dx
&\leq&
\int_{\mathbb R^n}\Big[\Big(\sum_j (\mathcal M(\chi_{B_j})(x))^r\Big)^{1/r}\Big]^{rq}\varphi(x,1)dx\\
&\leq&
C_{s,q}\int_{\mathbb R^n}\Big[\Big(\sum_j (\mathcal \chi_{B_j}(x))^r\Big)^{1/r}\Big]^{rq}\varphi(x,1)dx\\
&\leq&
C_{s,q} L\int_{\Omega} \varphi(x,1)dx\\
&\leq&
C \max\{1/\lambda, 1/\lambda^p\}\int_{\mathbb R^n}\varphi(x, f^*_m(x))dx
\end{eqnarray*}
for some $p\in(0,i(\varphi))$ since $\varphi\in \mathbb A_q\subset \mathbb A_{rq}$ and $f^*_m>\lambda$ on $\Omega$. Combining this, Lemma C and the H\"older inequality, we obtain
\begin{eqnarray*}
\int_{\mathbb R^n}[g^*_m(x)]^q \varphi(x,1)dx 
&\leq&
C\lambda^q \max\{1/\lambda, 1/\lambda^p\}\int_{\mathbb R^n}\varphi(x, f^*_m(x))dx+ C\int_{\Omega^c}[f^*_m(x)]^q \varphi(x,1)dx\\
&\leq&
C\lambda^q \max\{1/\lambda, 1/\lambda^p\}\int_{\mathbb R^n}\varphi(x, f^*_m(x))dx.
\end{eqnarray*}
 Here we used $f^*_m(x)\leq \lambda$ and $\varphi(x,\lambda)/\lambda^q\leq C \varphi(x,f^*_m(x))/[f^*_m(x)]^q$ for all $x\in \Omega^c$.
\end{proof}

\begin{Proposition}\label{the density of Lebesgue spaces in new Hardy spaces}
For any $q\in (q(\varphi),\infty)$ and $m\geq m(\varphi)$. The subspace $L^q_{\varphi(\cdot,1)}(\mathbb R^n)\cap H^\varphi_m(\mathbb R^n)$ is dense in $ H^\varphi_m(\mathbb R^n)$.
\end{Proposition}
\begin{proof}
 Let $f$ be an arbitrary element  in $H^\varphi_m(\mathbb R^n)$. For each $\lambda>0$, let $f=g^\lambda+\sum_j b^\lambda_j$ be the Calderon-Zygmund decomposition of $f$ of degree $m(\varphi)$, and height $\lambda$ associated with $f^*_m$. Then by Lemma \ref{basis lemma for atomic decomposition 1} and Lemma \ref{basis lemma for atomic decomposition 2}, $g^\lambda\in L^q_{\varphi(\cdot,1)}(\mathbb R^n)\cap H^\varphi_m(\mathbb R^n)$, moreover,
$$\int_{\mathbb R^n}\varphi(x, (g^\lambda- f)^*_m(x))dx\leq C \int_{f^*_m(x)>\lambda}\varphi(x, f^*_m(x))dx\to 0,$$
as $\lambda\to\infty$. Consequently, $\|g^\lambda-f\|_{ H^\varphi_m}\to 0$ as $\lambda\to\infty$ by Lemma \ref{the basis lemma 2}. Thus $L^q_{\varphi(\cdot,1)}(\mathbb R^n)\cap H^\varphi_m(\mathbb R^n)$ is dense in $H^\varphi_m(\mathbb R^n)$. 
\end{proof}

\subsection{The atomic decompositions $H^\varphi_m(\mathbb R^n)$}

 Recall that $m, s\geq m(\varphi)$, and $f$ is a distribution such that $f^*_m\in L^\varphi(\mathbb R^n)$. For each $k\in\mathbb Z$, let $f=g^k+\sum_j b_j^k$ be the Calder\'on-Zygmund decomposition of $f$ of degree $s$ and height $2^k$ associated with $f^*_m$. We shall label all the ingredients in this construction as in subsection \ref{subsection for Calderon-Zygmund decompositions}, but with superscript $k's$: for example,
$$\Omega^k=\{x\in\mathbb R^n: f^*_m(x)>2^k\},\qquad b_j^k=(f-P_j^k)\zeta_j^k,\qquad B_j^k=B(x_j^k,r_j^k).$$ 

Moreover, for each $k\in \mathbb Z$, and $i,j,$ let $P^{k+1}_{i,j}$ be the orthogonal projection of $(f-P_j^{k+1})\zeta_i^k$ onto $\mathcal P_s$ with respect to the norm associated to $\zeta_j^{k+1}$, namely, the unique element of $\mathcal P_s$ such that for all $Q\in \mathcal P_s$,
$$\int_{\mathbb R^n}(f(x)-P_j^{k+1}(x))\zeta_i^k(x) Q(x)\zeta_j^{k+1}(x)dx=\int_{\mathbb R^n}P^{k+1}_{i,j}(x)Q(x)\zeta_j^{k+1}(x)dx.$$

For convenience, we set $\hat B_j^k=B(x_j^k, 2r_j^k)$. Then we have the following lemma.

\begin{lemmad} [see \cite{FoS}, Chapter 3]

i) If $\hat B_j^{k+1}\cap \hat B_i^k\ne\emptyset$, then $r_j^{k+1}< 4r_i^k$ and $\hat B_j^{k+1}\subset B(x_i^k, 18r_i^k)$.\\
ii) For each $j$ there are at most $L$ (depending only on $n$ as in last section) values of $i$ such that $\hat B_j^{k+1}\cap \hat B_i^k\ne\emptyset$.\\
iii) There is a constant $C>0$, independent of $f,i,j$, and $k$, such that
$$\sup\limits_{x\in\mathbb R^n}|P^{k+1}_{i,j}(x)\zeta_j^{k+1}(x)|\leq C 2^{k+1}.$$
iv) For every $k\in \mathbb Z$, $\sum_i(\sum_jP^{k+1}_{i,j}\zeta^{k+1}_j)=0$, where the series converges pointwise and in $\mathcal S'(\mathbb R^n)$.
\end{lemmad}

We now give the necessary estimates for proving that $H^\varphi_m(\mathbb R^n)\subset H^{\varphi,\infty,s}_{\rm at}(\mathbb R^n)$, $m\geq s\geq m(\varphi)$, and the inclusion is continuous.

\begin{Lemma}\label{the fundamental lemma of Ky}
Let $f\in H^\varphi_m(\mathbb R^n)$, and for each $k\in\mathbb Z$, set 
$$\Omega^k=\{x\in\mathbb R^n: f^*_m(x)>2^k\}.$$
Then for any $\lambda>0$, there exists a constant $C$, independent of $f$ and $\lambda$, such that
$$\sum_{k=-\infty}^\infty \varphi\Big(\Omega^k, \frac{2^k}{\lambda}\Big)\leq C\int_{\mathbb R^n}\varphi\Big(x,\frac{f^*_m(x)}{\lambda}\Big)dx.$$
\end{Lemma}

\begin{proof}
Let $p\in (0,i(\varphi))$ and $C_p$ is such that (\ref{lower type}) holds.  We now set $N_0=[(\log_2 C_p)/p]+1$ so that $2^{N_0p}>C_p$.  For each $\ell\in\mathbb N,0\leq \ell\leq N_0-1$, we consider the sequence $U_m^\ell=\sum_{k=-m}^m\varphi\Big(\Omega^{N_0k+\ell}, \frac{2^{N_0k+\ell}}{\lambda}\Big)$. Obviously, $\{U_m^\ell\}_{m\in\mathbb N}$ is an increasing sequence. Moreover, for any $m\in\mathbb N$, 
\begin{eqnarray*}
U_m^\ell
&=&
\sum_{k=-m}^m\varphi\Big(\Omega^{N_0(k+1)+\ell}, \frac{2^{N_0k+\ell}}{\lambda}\Big)+\sum_{k=-m}^m \Big\{\varphi\Big(\Omega^{N_0k+\ell}, \frac{2^{N_0k+\ell}}{\lambda}\Big)-\varphi\Big(\Omega^{N_0(k+1)+ \ell}, \frac{2^{N_0k+\ell}}{\lambda}\Big)\Big\}\\
&\leq&
C_p\frac{1}{2^{N_0p}}\Big\{U_m^\ell +\varphi\Big(\Omega^{N_0(m+1)+\ell}, \frac{2^{N_0(m+1)+\ell}}{\lambda}\Big)+ \varphi\Big(\Omega^{N_0(-m)+\ell}, \frac{2^{N_0(-m)+\ell}}{\lambda}\Big)\Big\}+\\
&&
\qquad\qquad\qquad\qquad\qquad\qquad\qquad+\sum_{k=-m}^m\int_{\Omega^{N_0k+\ell}\setminus\Omega^{N_0(k+1)+\ell}}\varphi\Big(x, \frac{f^*_m(x)}{\lambda}\Big)dx\\
&\leq&
\frac{C_p}{2^{N_0p}}U_m^\ell+\Big(2\frac{C_p}{2^{N_0p}}+1\Big)\int_{\mathbb R^n}\varphi\Big(x, \frac{f^*_m(x)}{\lambda}\Big)dx.
\end{eqnarray*}
This implies that $U_m^\ell\leq \frac{3}{1-C_p/(2^{N_0p})}\int_{\mathbb R^n}\varphi\Big(x, \frac{f^*_m(x)}{\lambda}\Big)dx$. Consequently,
$$\sum_{k=-\infty}^\infty \varphi\Big(\Omega^k, \frac{2^k}{\lambda}\Big)=\sum_{\ell=0}^{N_0-1}\lim\limits_{m\to\infty}U_m^\ell\leq C \int_{\mathbb R^n}\varphi\Big(x, \frac{f^*_m(x)}{\lambda}\Big)dx,$$
where $C=\frac{3N_0}{1-C_p/(2^{N_0p})}$ independent of $f$ and $\lambda$.
\end{proof}

\begin{Theorem}\label{part 2 of atomic decomposition}
Let $m\geq s\geq m(\varphi)$. Then, $H^\varphi_m(\mathbb R^n)\subset  H^{\varphi,\infty,s}_{\rm at}(\mathbb R^n)$ and the inclusion is continuous.
\end{Theorem}
\begin{proof}
 Suppose first that $f\in L^q_{\varphi(\cdot,1)}(\mathbb R^n)\cap H^\varphi_m(\mathbb R^n)$ for some $q\in (q(\varphi),\infty)$. Let $f=g^k+\sum_j b_j^k$ be the Calder\'on-Zygmund decompositions of $f$ of degree $s$ with height $2^k$, for $k\in\mathbb Z$ associated with $f^*_m$. By Proposition \ref{the density of Lebesgue spaces in new Hardy spaces}, $g^k\to f$ in $H^\varphi_m(\mathbb R^n)$ as $k\to\infty$, while by \cite[Lemma 4.10]{BLYZ}, $g^k\to 0$ uniformly as $k\to -\infty$ since $f\in L^q_{\varphi(\cdot,1)}(\mathbb R^n)$. Therefore, $f=\sum_{-\infty}^\infty(g^{k+1}-g^k)$  in $\mathcal S'(\mathbb R^n)$. Using \cite[Lemma 3.27]{FoS} together with the equation $\sum_i \zeta_i^k b_j^{k+1}=\chi_{\Omega^k} b_j^{k+1}=b_j^{k+1}$ by supp$b_j^{k+1}\subset \Omega^{k+1}\subset \Omega^k$, we get
\begin{eqnarray*}
g^{k+1}-g^k
&=&
(f-\sum_j b_j^{k+1})-(f-\sum_i b_i^k)\\
&=&
\sum_i b_i^k-\sum_j b_j^{k+1}+\sum_i\sum_j P_{i,j}^{k+1}\zeta_j^{k+1}\\
&=&
\sum_i\Big[b_i^k-\sum_j\Big(\zeta_i^k b_j^{k+1}-P_{i,j}^{k+1}\zeta_j^{k+1}\Big)\Big]\\
&=&
\sum_i h^k_i
\end{eqnarray*}
where all the series converge  in $\mathcal S'(\mathbb R^n)$ and almost everywhere. Furthermore,
\begin{equation}\label{part 2 of atomic decomposition 1}
h^k_i=(f-P_i^k)\zeta_i^k- \sum_j \Big((f-P_j^{k+1})\zeta_i^k-P_{i,j}^{k+1}\Big)\zeta_j^{k+1}.
\end{equation}

From this formula it is obvious that $\int_{\mathbb R^n}h_i^k(x)P(x)dx=0$ for all $P\in\mathcal P_s$. Moreover, $h_i^k=\zeta_i^k f\chi_{(\Omega^{k+1})^c}-P_i^k \zeta_i^k+\zeta_i^k\sum_j P_j^{k+1}\zeta_j^{k+1}+\sum_jP_{i,j}^{k+1}\zeta_j^{k+1}$, by $\sum_j\zeta_j^{k+1}=\chi_{\Omega^{k+1}}$. But $|f(x)|\leq C(m)  f^*_m(x)\leq C 2^{k+1}$ for almost every $x\in (\Omega^{k+1})^c$, so by Lemmas 3.8 and 3.26 of \cite{FoS}, and $\sum_j\zeta_j^{k+1}\leq L$,
\begin{equation}\label{part 2 of atomic decomposition 2}
\|h^k_i\|_{L^\infty}\leq C 2^{k+1}+ C 2^k+ C L 2^{k+1}+ C L 2^{k+1}\leq C 2^k,
\end{equation}

 Lastly, since $P_{i,j}^{k+1}=0$ unless $\hat B_i^k\cap \hat B_j^{k+1}\ne\emptyset$, it follows from (\ref{part 2 of atomic decomposition 1}) and \cite[Lemma 3.24]{FoS}, that $h_i^k$ is supported in $B(x_i^k, 18r_i^k)$. Thus $h_i^k$ is a multiple of $(\varphi,\infty,s)$-atom. Moreover, by (\ref{part 2 of atomic decomposition 2}) and Lemma \ref{the fundamental lemma of Ky}, for any $\lambda>0$,
\begin{eqnarray*}
\sum_{k\in \mathbb Z}\sum_i \varphi\Big(B(x_i^k, 18 r_i^k),\frac{\|h_i^k\|_{L^\infty}}{\lambda}\Big)
&\leq&
\sum_{k\in\mathbb Z} L\varphi(\Omega^k, C 2^k/\lambda)\\
&\leq&
C\int_{\mathbb R^n}\varphi\Big(x, \frac{f^*_m(x)}{\lambda}\Big)dx<\infty.
\end{eqnarray*}

 Thus the series $\sum_{k\in\mathbb Z}\sum_i h_i^k$ converges in $ H^{\varphi,\infty,s}_{\rm at}(\mathbb R^n)$ and defines an atomic decomposition of $f$. Moreover,
\begin{eqnarray*}
\sum_{k\in\mathbb Z}\sum_i \varphi \Big(B(x_i^k, 18 r_i^k), \frac{\|h_i^k\|_{L^\infty}}{\|f\|_{ H^\varphi_m}}\Big)
&\leq&
C\int_{\mathbb R^n}\varphi\Big(x, \frac{f^*_m(x)}{\|f\|_{ H^\varphi_m}}\Big)dx\\
&\leq&
C.
\end{eqnarray*}
Consequently, $\|f\|_{H^{\varphi,\infty,s}_{\rm at}}\leq \Lambda_\infty(\{h_i^k\})\leq C\|f\|_{H^\varphi_m}$ by Lemma \ref{the basis lemma 3}.

Now, let $f$ be an arbitrary element of $H^\varphi_m(\mathbb R^n)$. There exists a sequence $\{f_\ell\}_{\ell\geq 1}\subset L^q_{\varphi(\cdot,1)}(\mathbb R^n)\cap H^\varphi_m(\mathbb R^n)$ such that $f=\sum_{\ell=1}^\infty f_\ell$ in $H^\varphi_m(\mathbb R^n)$ (thus in $\mathcal S'(\mathbb R^n)$) and $\|f_\ell\|_{ H^\varphi_m}\leq 2^{2-\ell}\|f\|_{H^\varphi_m}$ for any $\ell\geq 1$. For any $\ell\geq 1$, let $f_\ell=\sum_j b_{j,\ell}$ be the atomic decomposition of $f_\ell$, with supp $b_{j,\ell}\subset B_{j,\ell}$ constructed above. Then $f=\sum_{\ell=1}^\infty\sum_j b_{j,\ell}$ is an atomic decomposition of $f$, and
\begin{eqnarray*}
\sum_{\ell=1}^\infty\sum_j \varphi\Big(B_{j,\ell}, \frac{\|b_{j,\ell}\|_{L^\infty}}{\|f\|_{H^\varphi_m}}\Big)
&\leq&
\sum_{\ell=1}^\infty\sum_i \varphi\Big(B_{j,\ell}, \frac{\|b_{j,\ell}\|_{L^\infty}}{2^{\ell-2}\|f_\ell\|_{H^\varphi_m}}\Big)\\
&\leq&
\sum_{\ell=1}^\infty C_p\frac{1}{(2^{\ell-2})^p}=:C,
\end{eqnarray*}
where $C_p$ is such that (\ref{lower type}) holds. Thus $f\in H^{\varphi,\infty,s}_{\rm at}(\mathbb R^n)$, moreover,
$$\|f\|_{ H^{\varphi,\infty,s}_{\rm at}}\leq \Lambda_\infty(\{b_{j,\ell}\})\leq C \|f\|_{H^\varphi_m}$$
by Lemma \ref{the basis lemma 3}. This completes the proof.
\end{proof}

\begin{proof}[\bf Proof of Theorem \ref{atomic decomposition for new Hardy spaces}]
  By Theorem \ref{part 1 for the atomic decomposition} and Theorem \ref{part 2 of atomic decomposition}, we obtain
$$H^{\varphi,\infty,s}_{\rm at}(\mathbb R^n)\subset H^{\varphi,q,s}_{\rm at}(\mathbb R^n)\subset H^{\varphi,q,m(\varphi)}_{\rm at}(\mathbb R^n)\subset H^\varphi(\mathbb R^n)\subset H^\varphi_s(\mathbb R^n)\subset H^{\varphi,\infty,s}_{\rm at}(\mathbb R^n)$$
and the inclusions are continuous. Thus $ H^\varphi(\mathbb R^n)= H^{\varphi,q,s}_{\rm at}(\mathbb R^n)$ with equivalent norms.  
\end{proof}

\section{Dual spaces}

In this section, we give the proof of Theorem \ref{the dual theorem for new Hardy spaces}. In order to do this, we need the below  lemma, which can be seen as a consequence of the fact that $\varphi(\cdot,t)$ is {\sl uniformly locally integrable}. We omit the details here.

\begin{Lemma}\label{a lemma for the dual theorem for new Hardy spaces}
Given a ball $B$, and $\{B_j\}_j$ be a sequence of measurable subsets of $B$ such that $\lim\limits_{j\to\infty}|B_j|=0$. Then the following holds
$$\lim\limits_{j\to\infty}\sup\limits_{t>0}\frac{\varphi(B_j,t)}{\varphi(B,t)}=0.$$
\end{Lemma}

We next note that if $\mathfrak b\in BMO^\varphi(\mathbb R^n)$ is real-valued and 
\begin{equation*}
\mathfrak b_N(x)=
\begin{cases}
N &  \quad{\rm if}\quad \mathfrak b(x)> N,\\
\mathfrak b(x) & \quad{\rm if}\quad |\mathfrak b(x)|\leq N,\\
-N & \quad{\rm if}\quad\mathfrak b(x)< -N,
\end{cases}
\end{equation*}
then by using the fact
$$\|f\|_{BMO^\varphi}\leq \sup\limits_{B\rm -ball}\frac{1}{\|\chi_B\|_{L^\varphi}}\frac{1}{|B|}\int_{B}\int_{B} |f(x)-f(y)|dxdy\leq 2\|f\|_{BMO^\varphi},$$
we obtain that  $\|\mathfrak b_N\|_{BMO^\varphi}\leq 2 \|\mathfrak b\|_{BMO^\varphi}$ for all $N>0$.

\begin{proof}[\bf Proof of Theorem \ref{the dual theorem for new Hardy spaces}]
  i) It is sufficient to prove it for $\mathfrak b\in BMO^\varphi(\mathbb R^n)$ real-valued since $\mathfrak b\in BMO^\varphi(\mathbb R^n)$ iff $\mathfrak b= \mathfrak b_1+ i \mathfrak b_2$ with $\mathfrak b_j\in BMO^\varphi(\mathbb R^n)$ real-valued, $j=1,2$, moreover 
$$\|\mathfrak b\|_{BMO^\varphi}\approx \|\mathfrak b_1\|_{BMO^\varphi}+ \|\mathfrak b_2\|_{BMO^\varphi}.$$

Suppose first that $\mathfrak b\in BMO^\varphi(\mathbb R^n)\cap L^\infty(\mathbb R^n)$. Then, the functional 
$$L_{\mathfrak b}(f)= \int_{\mathbb R^n} f(x)\mathfrak b(x)dx$$
is well defined for any $f\in L_0^\infty(\mathbb R^n)$ since $\mathfrak b\in L^1_{\rm loc}(\mathbb R^n)$. 

Furthermore, since $f\in L_0^\infty(\mathbb R^n)\subset L^2(\mathbb R^n)\cap H^1(\mathbb R^n)$, we remark that the atomic decomposition $f=\sum_{k\in\mathbb Z}\sum_i h^k_i$ in the proof of Theorem \ref{part 2 of atomic decomposition} is also the classical atomic decomposition of $f$ in $H^1(\mathbb R^n)$, so that the series converge in $H^1(\mathbb R^n)$ and thus in $L^1(\mathbb R^n)$. Combining this with the fact $\mathfrak b\in L^\infty(\mathbb R^n)$, we obtain
$$L_{\mathfrak b}(f)= \int_{\mathbb R^n} f(x)\mathfrak b(x)dx=\sum_{k\in\mathbb Z}\sum_i \int_{\mathbb R^n} h_i^k(x) \mathfrak b(x)dx.$$

Therefore, by Lemma \ref{the basis lemma 4} and the proof of Theorem \ref{part 2 of atomic decomposition},
\begin{eqnarray*}
|L_{\mathfrak b}(f)|= \Big|\int_{\mathbb R^n} f(x)\mathfrak b(x)dx\Big|
&\leq& \sum_{k\in\mathbb Z} \sum_i \Big|\int_{\mathbb R^n} h_i^k(x) \mathfrak b(x)dx\Big|\\
&=&  \sum_{k\in\mathbb Z} \sum_i \Big|\int_{B(x_i^k, 18 r_i^k)} h_i^k(x) (\mathfrak b(x)- \mathfrak b_{B(x_i^k, 18 r_i^k)}(x))dx\Big|\\
&\leq&
\|\mathfrak b\|_{BMO^\varphi} \sum_{k\in\mathbb Z} \sum_i \|h_i^k\|_{L^\infty}\|\chi_{B(x_i^k, 18 r_i^k)}\|_{L^\varphi}\\
&\leq&
C \|\mathfrak b\|_{BMO^\varphi} \Lambda_\infty(\{h_i^k\})\\
&\leq&
C \|\mathfrak b\|_{BMO^\varphi} \|f\|_{H^\varphi}.
\end{eqnarray*}

Now, let $\mathfrak b$ be an arbitrary element in $BMO^\varphi(\mathbb R^n)$. For any $f\in L^\infty_0(\mathbb R^n)$, it is clear that  $|f \mathfrak b_\ell| \leq |f \mathfrak b|\in L^1(\mathbb R^n)$ for every $\ell\geq 1$, and $f(x)\mathfrak b_\ell(x)\to f(x)\mathfrak b(x)$, as $\ell\to\infty$, for almost every $x\in\mathbb R^n$. Therefore, by the dominated convergence theorem of Lebesgue, we obtain
$$|L_{\mathfrak b}(f)|= \Big|\int_{\mathbb R^n} f(x)\mathfrak b(x)dx\Big|=\lim\limits_{\ell\to\infty}\Big|\int_{\mathbb R^n} f(x)\mathfrak b_\ell(x)dx\Big|\leq C  \|\mathfrak b\|_{BMO^\varphi} \|f\|_{H^\varphi},$$
since $\|\mathfrak b_\ell\|_{BMO^\varphi}\leq 2 \|\mathfrak b\|_{BMO^\varphi}$ for all $\ell\geq 1$.

Because of the density of $L^\infty_0(\mathbb R^n)$ in $H^\varphi(\mathbb R^n)$, the functional $L_\mathfrak b$ can be extended to a bounded functional on $H^\varphi(\mathbb R^n)$, moreover, $\|L_\mathfrak b\|_{(H^\varphi)^*}\leq C\|\mathfrak b\|_{BMO^\varphi}$.

ii) Conversely, suppose $L$ is a continuous linear functional on $H^\varphi(\mathbb R^n)\equiv H^{\varphi,q, 0}_{\rm at}(\mathbb R^n)$ for some $q\in (q(\varphi), \infty)$. For any ball $B$, denote by $ L^q_{\varphi,0}(B)$ the subspace of $L^q_\varphi(B)$ defined by
$$ L^q_{\varphi,0}(B):=\Big\{f\in L^q_\varphi(B): \int_{\mathbb R^n}f(x)dx=0\Big\}.$$
Obviously, if $B_1\subset B_2$ then
\begin{equation}\label{the dual theorem for new Hardy spaces 1}
L^q_\varphi(B_1)\subset L^q_\varphi(B_2)\quad\mbox{and}\quad  L^q_{\varphi,0}(B_1)\subset L^q_{\varphi,0}(B_2).
\end{equation}
Moreover, when $f\in  L^q_{\varphi,0}(B)\setminus \{0\}$, $a(x)= \|\chi_B\|_{L^\varphi}^{-1}\|f\|^{-1}_{L^q_\varphi(B)}f(x)$ is a $(\varphi,q, 0)$-atom, thus $f\in H^{\varphi,q, 0}_{\rm at}(\mathbb R^n)$ and 
$$\|f\|_{ H^{\varphi,q, 0}_{\rm at}}\leq \|\chi_B\|_{L^\varphi}\|f\|_{L^q_\varphi(B)}.$$
Since $L\in (H^{\varphi,q, 0}_{\rm at}(\mathbb R^n))^*$, by the above,
$$|L(f)|\leq \|L\|_{( H^{\varphi,q, 0}_{\rm at})^*}\|f\|_{ H^{\varphi,q, 0}_{\rm at}}\leq \|L\|_{(H^{\varphi,q,0}_{\rm at})^*}\|\chi_B\|_{L^\varphi}\|f\|_{L^q_\varphi(B)},$$
for all $f\in L^q_{\varphi,0}(B)$. Therefore, $L$ provides a bounded linear functional on $ L^q_{\varphi,0}(B)$
which can be extended by the Hahn-Banach theorem to the whole space $L^q_\varphi(B)$ without increasing its norm. On the other hand, by Lemma \ref{a lemma for the dual theorem for new Hardy spaces} and Lebesgue-Nikodym theorem, there exists $h\in L^1(B)$ such that
$$L(f)=\int_{\mathbb R^n} f(x) h(x)dx,$$
for all $f\in L^\infty_{\varphi,0}(B)$. 

We now take a sequence of balls $\{B_j\}_{j\geq 1}$ such that $B_1\subset B_2\subset\cdots\subset B_j\subset\cdots$ and $\cup_j B_j=\mathbb R^n$. Then, there exists a sequence $\{h_j\}_{j\geq 1}$ such that
$$h_j\in L^1(B_j)\quad\mbox{and}\quad L(f)=\int_{\mathbb R^n} f(x)h_j(x)dx,$$
for all $f\in  L^\infty_{\varphi,0}(B_j), j=1,2,...$ Hence, for all $f\in  L^\infty_{\varphi,0}(B_1)\subset  L^\infty_{\varphi,0}(B_2)$ (by (\ref{the dual theorem for new Hardy spaces 1})),
$$\int_{\mathbb R^n} f(x)(h_1(x)-h_2(x))dx=\int_{\mathbb R^n} f(x)h_1(x)dx - \int_{\mathbb R^n} f(x)h_2(x)dx= L(f)-L(f)=0.$$
As $f_{B_1}=0$ if $f\in  L^\infty_{\varphi,0}(B_1)$, we have
$$\int_{\mathbb R^n} f(x)\Big((h_1(x)-h_2(x))- (h_1-h_2)_{B_1}\Big)dx=0$$
for all $f\in  L^\infty_{\varphi,0}(B_1)$, and thus for $f\in L^\infty_\varphi(B_1)$. Hence,
$$h_1(x)-h_2(x)=(h_1-h_2)_{B_1}\;,\; a.e\; x\in B_1.$$

By the similar arguments, we also obtain
\begin{equation}\label{the dual theorem for new Hardy spaces 2}
h_j(x)-h_{j+1}(x)= (h_j-h_{j+1})_{B_j}
\end{equation}
a.e $x\in B_j, j=2,3,...$ Consequently, if we define the sequence $\{\widetilde h_j\}_{j\geq 1}$ by
\begin{equation*}
\begin{cases}
\widetilde h_1=h_1\\
\widetilde h_{j+1}=h_{j+1}+ (\widetilde h_j-h_{j+1})_{B_j}&, \quad j= 1,2,...
\end{cases}
\end{equation*}
then it follows from (\ref{the dual theorem for new Hardy spaces 2}) that
$$\widetilde h_j\in L^1(B_j)\quad\mbox{and}\quad \widetilde h_{j+1}(x)=\widetilde h_j(x)$$
a.e $x\in B_j, j=1,2,...$ Thus, we can define the function $\mathfrak b$ on $\mathbb R^n$ by
$$\mathfrak b(x)=\widetilde h_j(x)$$
if $x\in B_j$ for some $j\geq 1$ since $B_1\subset B_2\subset\cdots\subset B_j\subset\cdots$ and $\cup_j B_j=\mathbb R^n$.

Let us now show that $\mathfrak b\in BMO^\varphi(\mathbb R^n)$ and
$$L(f)=\int_{\mathbb R^n} f(x)\mathfrak b(x)dx,$$
for all $f\in L^\infty_0(\mathbb R^n)$.

 Indeed, for any $f\in L^\infty_0(\mathbb R^n)$, there exists $j\geq 1$ such that $f\in L^\infty_{\varphi,0}(B_j)$. Hence,
$$L(f)=\int_{\mathbb R^n} f(x)\widetilde h_j(x)dx= \int_{ B_j} f(x)\widetilde h_j(x)dx=\int_{\mathbb R^n} f(x)\mathfrak b(x)dx.$$

On the other hand, for all ball $B$, one consider $f={\rm sign}(\mathfrak b- \mathfrak b_B)$ where sign$\xi=\overline\xi/|\xi|$ if $\xi\ne 0$ and sign$0=0$. Then, 
$$a=\frac{1}{2}\|\chi_B\|_{L^\varphi}^{-1}(f- f_B)\chi_B$$
is a $(\varphi, \infty, 0)$-atom. Consequently, 
\begin{eqnarray*}
|L(a)|
&=&
\frac{1}{2}\|\chi_B\|_{L^\varphi}^{-1}\Big|\int_{\mathbb R^n}\mathfrak b(x)(f(x)- f_B)\chi_B(x)dx\Big|\\
&=&
\frac{1}{2}\frac{1}{\|\chi_B\|_{L^\varphi}}\Big|\int_B (\mathfrak b(x)- \mathfrak b_B)f(x)dx\Big|\\
&=&
\frac{1}{2}\frac{1}{\|\chi_B\|_{L^\varphi}}\int_B |\mathfrak b(x)- \mathfrak b_B|dx\\
&\leq&
\|L\|_{(H^\varphi)^*}\|a\|_{H^\varphi}\leq C \|L\|_{( H^\varphi)^*}
\end{eqnarray*}
since $L\in (H^\varphi(\mathbb R^n))^*$ and Corollary 5.2. As $B$ is arbitrary, the above implies $\mathfrak b\in BMO^\varphi(\mathbb R^n)$ and
$$\|\mathfrak b\|_{BMO^\varphi}\leq C \|L\|_{(H^\varphi)^*}.$$

The uniqueness (in the sense $\mathfrak b= \widetilde{\mathfrak b}$ if $\mathfrak b- \widetilde{\mathfrak b}= \rm const$) of the function $\mathfrak b$ is clear. And thus the proof is finished.
\end{proof}

\section{The class of pointwise multipliers for $BMO(\mathbb R^n)$}

 In this subsection, we give as an interesting application that the class of pointwise multipliers for $BMO(\mathbb R^n)$ is just the dual of $L^1(\mathbb R^n)+ H^{\rm log}(\mathbb R^n)$ where $ H^{\rm log}(\mathbb R^n)$ is a Hardy space of Musielak-Orlicz type related to the Musielak-Orlicz function $\theta(x,t)=\frac{t}{\log(e+|x|)+ \log(e+t)}$.

We first introduce {\sl log-atoms}. A measurable function $a$ is said to be {\sl log-atom} if it satisfies the following three conditions

$\bullet$ $a$ supported in $B$ for some ball $B$ in $\mathbb R^n$,

$\bullet$ $\|a\|_{L^\infty}\leq \displaystyle\frac{\log (e+\frac{1}{|B|})+ \sup_{x\in B}\log(e+|x|)}{|B|}$,

$\bullet$ $\int_{\mathbb R^n}a(x)dx=0$.\\

To prove Theorem \ref{a theorem on the class of pointwise multipliers for BMO}, we need the following two propositions.

\begin{Proposition}\label{equivalence of log-atoms}
There exists a positive constant $C$ such that if $f$ is a $\theta$-atom (resp., {\sl log-atom}) then $C^{-1}f$ is a {\sl log-atom} (resp., $\theta$-atom). 
\end{Proposition}

\begin{Proposition}\label{equivalence of BMO type spaces}
On  $BMO^{\rm log}(\mathbb R^n)$, we have
$$\|f\|_{BMO^{\rm log}}\approx \sup\limits_{B\rm -ball}\frac{\log (e+\frac{1}{|B|})+ \sup_{x\in B}\log(e+|x|)}{|B|}\int_{B}|f(x)-f_B|dx<\infty.$$
\end{Proposition}

We first note that $\theta$ is a {\sl growth function} that satisfies $n q(\theta)< (n+1) i(\theta)$ in Theorem \ref{the dual theorem for new Hardy spaces}. More precisely, $\theta\in \mathbb A_1$ and $\theta(x,\cdot)$ is concave with $i(\theta)=1$.

\begin{proof}[\bf Proof of Proposition \ref{equivalence of log-atoms}]
Let $f$ be a {\sl log-atom}. By the above remark, to prove that there exists a constant $C>0$ (independent of $f$ and which may change from line to line) such that $C^{-1}f$ is a $\theta$-atom, it is sufficient to show that there exists a constant $C>0$ such that
$$\int_B\theta\Big(x, \frac{\log (e+\frac{1}{|B|})+ \sup_{x\in B}\log(e+|x|)}{|B|}\Big)dx\leq C$$
or, equivalently, 
$$\frac{\frac{\log (e+\frac{1}{|B|})+ \sup_{x\in B}\log(e+|x|)}{|B|}}{\log(e+ \frac{\log (e+\frac{1}{|B|})+ \sup_{x\in B}\log(e+|x|)}{|B|})+ \sup_{x\in B}\log(e+|x|)}|B|\leq C,$$
since $\theta\in \mathbb A_1$. However, the last inequality is obvious. 

Conversely, suppose that $f$ is a $\theta$-atom. Similarly, we need to  show that there exists a constant $C>0$ such that
$$\int_B\theta\Big(x, C\frac{\log (e+\frac{1}{|B|})+ \sup_{x\in B}\log(e+|x|)}{|B|}\Big)dx\geq 1$$
or, equivalently, 
$$\frac{C\frac{\log (e+\frac{1}{|B|})+ \sup_{x\in B}\log(e+|x|)}{|B|}}{\log(e+ C\frac{\log (e+\frac{1}{|B|})+ \sup_{x\in B}\log(e+|x|)}{|B|})+ \sup_{x\in B}\log(e+|x|)}|B|\geq 1.$$
However it is true. For instance we may take $C=3$. 
\end{proof}

\begin{proof}[\bf Proof of Proposition \ref{equivalence of BMO type spaces}]
 It is sufficient to show that there exists a constant $C>0$ such that
$$C^{-1} (|\log r|+ \log (e+ |x|))\leq \log\Big(e+ \frac{1}{|B(x,r)|}\Big)+ \sup\limits_{y\in B(x,r)}\log(e+|y|)\leq C ( |\log r|+ \log (e+ |x|)).$$

 The first inequality is easy and shall be omited. For the second, one first consider the 1 dimensional case. Then by symmetry, we just need  to prove that
$$\log(e+ 1/(b-a))+ \sup\limits_{x\in [a,b]}\log(e+ |x|)\leq C (|\log (b-a)/2|+ \log(e+ |a+b|/2))$$
for all $b>0, a\in [-b,b)\subset \mathbb R$. However, this follows from the basic two inequalities:
$$\log(e+ 1/(b-a))\leq 2(|\log (b-a)/2|+ \log(e+ |a+b|/2))$$
and
$$\log(e+b)\leq 5 \log(e+b)/2\leq 5(|\log (b-a)/2|+ \log(e+ |a+b|/2)).$$

For the general case $\mathbb R^n$, by the 1-dimensional result, we obtain
\begin{eqnarray*}
\log\Big(e+ \frac{1}{|B(x,r)|}\Big)&\leq& \frac{2^n}{c_n}\sum_{i=1}^n\log\Big(e+ \frac{1}{|[x_i-r, x_i+r]|}\Big)\\
&\leq&
 C \sum_{i=1}^n( |\log r|+ \log (e+ |x_i|))\\
&\leq&
 C ( |\log r|+ \log (e+ |x|))
\end{eqnarray*}
where $c_n=|B(0,1)|$, and
\begin{eqnarray*}
\sup\limits_{y\in B(x,r)}\log(e+|y|)&\leq& \sum_{i=1}^n \sup\limits_{y_i\in [x_i-r, x_i+r]}\log(e+|y_i|)\\
&\leq&
C \sum_{i=1}^n ( |\log r|+ \log (e+ |x_i|))\\
&\leq&
 C( |\log r|+ \log (e+ |x|))
\end{eqnarray*}
where $x=(x_1,...,x_n), y=(y_1,...,y_n)\in\mathbb R^n$. This finishes the proof.  
\end{proof}

\begin{proof}[\bf Proof of Theorem \ref{a theorem on the class of pointwise multipliers for BMO}]
 By Theorem \ref{atomic decomposition for new Hardy spaces}, Theorem \ref{the dual theorem for new Hardy spaces}, Proposition \ref{equivalence of log-atoms}, and Proposition \ref{equivalence of BMO type spaces}, we obtain $(H^{\rm log}(\mathbb R^n))^*\equiv BMO^{\rm log}(\mathbb R^n)$. We deduce that, the class of pointwise multipliers for $BMO(\mathbb R^n)$ is the dual of $L^1(\mathbb R^n)+ H^{\rm log}(\mathbb R^n)$.  
\end{proof}

\section{Finite atomic decompositions and their applications}

We first prove the finite atomic decomposition theorem.

\begin{proof}[\bf Proof of Theorem \ref{finite decomposition for new Hardy spaces}]
Obviously, $H^{\varphi,q,s}_{\rm fin}(\mathbb R^n)\subset H^\varphi(\mathbb R^n)$ and for all $f\in H^{\varphi,q,s}_{\rm fin}(\mathbb R^n)$,
$$\|f\|_{H^\varphi}\leq C \|f\|_{H^{\varphi,q,s}_{\rm fin}}.$$
Thus, we have to show that for every $q\in (q(\varphi),\infty)$ there exists a constant $C>0$ such that
$$\|f\|_{H^{\varphi,q,s}_{\rm fin}}\leq C \|f\|_{H^\varphi}$$
for all $f\in H^{\varphi,q,s}_{\rm fin}(\mathbb R^n)$ and that a similar estimate holds for $q=\infty$ and all $f\in H^{\varphi,\infty,s}_{\rm fin}(\mathbb R^n)\cap C(\mathbb R^n)$.

Assume that $q\in (q(\varphi),\infty]$, and by homogeneity, $f\in H^{\varphi,q,s}_{\rm fin}(\mathbb R^n)$ with $\|f\|_{H^\varphi}=1$. Notice that $f$ has compact support. Suppose that supp $f\subset B=B(x_0,r)$ for some ball $B$. Recall that, for each $k\in\mathbb Z$,
$$\Omega_k=\{x\in \mathbb R^n: f^*(x)> 2^k\}.$$

Clearly, $f\in L^{\overline q}_{\varphi(\cdot,1)}(\mathbb R^n)\cap H^\varphi(\mathbb R^n)$ where $\overline q=q$ if $q<\infty$, $\overline q= q(\varphi)+1$ if $q=\infty$. Hence, by Theorem \ref{part 2 of atomic decomposition}, there exists an atomic decomposition $f=\sum_{k\in\mathbb Z}\sum_i h_i^k\in H^{\varphi,\infty,s}_{\rm at}(\mathbb R^n)\subset H^{\varphi,q,s}_{\rm at}(\mathbb R^n)$ where the series converges in $\mathcal S'(\mathbb R^n)$ and almost everywhere. Moreover,
\begin{equation}\label{finite decomposition for new Hardy spaces 1}
\Lambda_q(\{h^k_i\})\leq \Lambda_\infty (\{h^k_i\})\leq C \|f\|_{H^\varphi}=C.
\end{equation}

On the other hand, it follows from the second step in the proof of Theorem 6.2 of \cite{BLYZ}  that there exists a constant $\widetilde C>0$, depending only on $m(\varphi)$, such that $f^*(x)\leq \widetilde C \inf\limits_{y\in B}f^*(y)$ for all $x\in B(x_0, 2r)^c$. Hence, we have
$$f^*(x)\leq \widetilde C \inf\limits_{y\in B}f^*(y)\leq \widetilde C \|\chi_B\|_{L^\varphi}^{-1}\|f^*\|_{L^\varphi}\leq \widetilde C \|\chi_B\|_{L^\varphi}^{-1}$$
for all $x\in B(x_0, 2r)^c$.  We now denote by $k'$ the largest integer $k$ such that $2^k< \widetilde C \|\chi_B\|_{L^\varphi}^{-1}$. Then,
\begin{equation}\label{finite decomposition for new Hardy spaces 2}
\Omega_k\subset B(x_0, 2r)\quad\mbox{for all}\; k>k'.
\end{equation}

Next we define the functions $g$ and $\ell$ by
$$g=\sum_{k\leq k'}\sum_i h^k_i\quad\mbox{and}\quad \ell =\sum_{k> k'}\sum_i h^k_i,$$
where the series converge in $\mathcal S'(\mathbb R^n)$ and almost everywhere. Clearly, $f= g+\ell$ and supp $\ell\subset \cup_{k>k'}\Omega_k\subset B(x_0, 2r)$ by (\ref{finite decomposition for new Hardy spaces 2}). Therefore, $g=f=0$ in $B(x_0, 2r)^c$, and thus supp $g\subset B(x_0,2r)$.

Let $1< \widetilde q< \frac{q}{q(\varphi)}$, then $\varphi\in \mathbb A_{q/\widetilde q}$. Consequently, 
$$\Big(\frac{1}{|B|}\int_B |f(x)|^{\widetilde q}dx\Big)^{1/\widetilde q}\leq C\left(\frac{1}{\varphi(B,1)}\int_B |f(x)|^q \varphi(x,1)dx\right)^{1/q}<\infty$$
by Lemma \ref{the basis lemma 5} if $q<\infty$ and it is trivial if $q=\infty$. Observe that supp $f\subset B$ and that $f$ has vanishing moments up to order $s$. By the above, we obtain that $f$ is a multiple of a classical $(1,\widetilde q, 0)$-atom and thus $f^*\in L^1(\mathbb R^n)$. Hence, it follows from (\ref{finite decomposition for new Hardy spaces 2}) that 
$$\int_{\mathbb R^n}\sum_{k>k'}\sum_i |h^k_i(x)x^\alpha|dx\leq C(|x_0|+ 2r)^s \sum_{k>k'} 2^k|\Omega_k|\leq C(|x_0|+ 2r)^s \|f^*\|_{L^1}<\infty,$$
for all $|\alpha|\leq s$. This together with the vanishing moments of $h^k_i$ implies that $\ell$ has vanishing moments up to order $s$ and thus so does $g$ by $g=f-\ell$.

In order to estimate the size of $g$ in $B(x_0, 2r)$, we recall that 
\begin{equation}\label{finite decomposition for new Hardy spaces 3}
\|h^k_i\|_{L^\infty}\leq C 2^k\;, \; \mbox{supp}\; h^k_i\subset B(x_i^k, 18r^k_i) \;\mbox{and}\; \sum_i \chi_{B(x_i^k, 18r^k_i)}\leq C. 
\end{equation}
Combining the above and the fact $\|\chi_B\|_{L^\varphi}\approx \|\chi_{B(x_0,2r)}\|_{L^\varphi}$, we obtain
$$\|g\|_{L^\infty}\leq C \sum_{k\leq k'}2^k\leq C 2^{k'}\leq C\widetilde C\|\chi_B\|_{L^\varphi}^{-1}\leq C \|\chi_{B(x_0,2r)}\|_{L^\varphi}^{-1}.$$
This proves that (see Definition \ref{atom})
\begin{equation}\label{finite decomposition for new Hardy spaces 4}
C^{-1}g \;\mbox{ is a}\; (\varphi,\infty,s){\rm -atom}.
\end{equation}

Now, we assume that $q\in (q(\varphi),\infty)$ and conclude the proof of (i). We first verify $\sum_{k>k'}\sum_i h^k_i\in L^q_\varphi(B(x_0,2r))$. For any $x\in \mathbb R^n$, since $\mathbb R^n=\cup_{k\in\mathbb Z}(\Omega_k\setminus \Omega_{k+1})$, there exists $j\in\mathbb Z$ such that $x\in \Omega_j\setminus \Omega_{j+1}$. Since supp $h^k_i\subset \Omega_k\subset \Omega_{j+1}$ for $k\geq j+1$,  it follows from (\ref{finite decomposition for new Hardy spaces 3}) that
$$\sum_{k>k'}\sum_i |h^k_i(x)|\leq C \sum_{k\leq j}2^k \leq C 2^j\leq C f^*(x).$$

Since $f\in L^q_\varphi(B)\subset L^q_\varphi(B(x_0,2r))$, we have $f^*\in  L^q_\varphi(B(x_0,2r))$. As $\varphi$ satisfies uniformly locally dominated convergence condition, we further obtain $\sum_{k>k'}\sum_{i} h^k_i$ converges to $\ell$ in $L^q_\varphi(B(x_0,2r))$.

Now, for any positive integer $K$, set $F_K=\{(i,k): k>k', |i|+ |k|\leq K\}$ and $\ell_K=\sum_{(i,k)\in F_K} h^k_i$. Observe that since $\sum_{k>k'}\sum_{i} h^k_i$ converges to $\ell$ in $L^q_\varphi(B(x_0,2r))$,  for any $\varepsilon>0$, if $K$ is large enough, we have $\varepsilon^{-1}(\ell-\ell_K)$ is a $(\varphi,q,s)$-atom. Thus, $f=g+ \ell_K + (\ell- \ell_K)$ is a finite linear atom combination of $f$. Then, it follows from (\ref{finite decomposition for new Hardy spaces 1}) and $(\ref{finite decomposition for new Hardy spaces 4})$ that
$$\|f\|_{H^{\varphi,q,s}_{\rm fin}}\leq C( C + \Lambda_q(\{h^k_i\}_{(i,k)\in F_K})+ \varepsilon)\leq C,$$
which ends the proof of (i).

To prove (ii), assume that $f$ is a continuous function in $H^{\varphi,\infty,s}_{\rm fin}(\mathbb R^n)$, and thus $f$ is uniformly continuous. Then, $h^k_i$ is continuous by examining its definition. Since $f$ is bounded, there exists a positive integer $k''>k'$ such that $\Omega_k=\emptyset$ for all $k>k''$. Consequently, $\ell= \sum_{k'<k\leq k''}\sum_{i} h^k_i$.

Let $\varepsilon>0$. Since $f$ is uniformly continuous, there exists  $\delta>0$ such that if $|x-y|<\delta$, then $|f(x)-f(y)|<\varepsilon$. Write $\ell= \ell_1^\varepsilon+ \ell_2^\varepsilon$ with 
$$\ell_1^\varepsilon\equiv \sum_{(i,k)\in F_1}h^k_i\quad \mbox{and}\quad \ell_2^\varepsilon\equiv \sum_{(i,k)\in F_2}h^k_i$$ where $F_1= \{(i,k): C r^k_i\geq \delta, k'< k\leq k''\}$ and $F_2= \{(i,k): C r^k_i< \delta, k'< k\leq k''\}$ with $C> 36$ the geometric constant (see \cite{MSV2}). Notice that the remaining part $\ell_1^\varepsilon$ will then be a finite sum. Since the atoms are continuous, $\ell_1^\varepsilon$ will be a continuous function. Furthermore, $\|\ell_2^\varepsilon\|_{L^\infty}\leq C (k''-k')\varepsilon$ (see also \cite{MSV2}). This means that one can write $\ell$ as the sum of one continuous term and of one which is uniformly arbitrarily small. Hence, $\ell$ is continuous, and so is $g = f -\ell$.

To find a finite atomic decomposition of $f$, we use again the splitting $\ell= \ell_1^\varepsilon+ \ell_2^\varepsilon$. By (\ref{finite decomposition for new Hardy spaces 1}), the part $\ell_1^\varepsilon$ is a finite sum of multiples of $(\varphi,\infty,s)$-atoms, and
\begin{equation}\label{finite decomposition for new Hardy spaces 5}
\|\ell_1^\varepsilon\|_{H^{\varphi,\infty,s}_{\rm fin}}\leq \Lambda_\infty (\{h^k_i\})\leq C \|f\|_{H^\varphi}=C.
\end{equation}

By $\ell, \ell_1^\varepsilon$ are continuous and have vanishing moments up to order $s$, and thus so does $\ell_2^\varepsilon= \ell- \ell_1^\varepsilon$. Moreover,  supp $\ell_2^\varepsilon\subset B(x_0, 2r)$ and $\|\ell_2^\varepsilon\|_{L^\infty}\leq C (k''-k')\varepsilon$. So we can choose $\varepsilon$ small enough such that $\ell_2^\varepsilon$ into an arbitrarily small multiple of a continuous $(\varphi,\infty,s)$-atom. Therefore, $f=g+ \ell_1^\varepsilon+ \ell_2^\varepsilon$ is a finite linear continuous atom combination of $f$. Then, it follows from (\ref{finite decomposition for new Hardy spaces 4}) and (\ref{finite decomposition for new Hardy spaces 5})  that
$$\|f\|_{H^{\varphi,\infty,s}_{\rm fin}}\leq C(\|g\|_{H^{\varphi,\infty,s}_{\rm fin}}+ \|\ell_1^\varepsilon\|_{H^{\varphi,\infty,s}_{\rm fin}}+ \|\ell_2^\varepsilon\|_{H^{\varphi,\infty,s}_{\rm fin}})\leq C.$$
This finishes the proof of (ii) and hence, the proof of Theorem 3.4. 

\end{proof}

Next we give the proof for Theorem \ref{the boundedness of sublinear operators on new Hardy spaces}.

\begin{proof}[\bf Proof of Theorem \ref{the boundedness of sublinear operators on new Hardy spaces}]
Suppose that the assumption (i) holds. For any $f\in H^{\varphi,q,s}_{\rm fin}(\mathbb R^n)$, by Theorem \ref{finite decomposition for new Hardy spaces}, there exists a finite atomic decomposition $f=\sum_{j=1}^k \lambda_j a_j$, where $a_j$'s are $(\varphi,q,s)$-atoms with supported in balls $B_j$'s, such that 
$$\Lambda_q(\{\lambda_j a_j\}_{j=1}^k)=\inf\left\{\lambda>0: \sum_{j=1}^k \varphi\Big(B_j,\frac{|\lambda_j| \|\chi_{B_j}\|_{L^\varphi}^{-1}}{\lambda}\Big)\leq 1\right\}\leq C\|f\|_{H^\varphi}.$$

Recall that, since $\varphi$ is of uniformly upper type $\gamma$, there exists a constant $C_\gamma>0$ such that
\begin{equation}\label{the boundedness of sublinear operators on new Hardy spaces 1}
\varphi(x, st)\leq C_\gamma s^\gamma \varphi(x, t)\;\;\mbox{for all}\; x\in\mathbb R^n, s\geq 1, t\in [0,\infty).
\end{equation}
If there exist $j_0\in \{1,...,k\}$ such that $C_\gamma|\lambda_{j_0}|^\gamma\geq \sum_{j=1}^k |\lambda_j|^\gamma$, then
$$\sum_{j=1}^k \varphi\left(B_j, \frac{|\lambda_j| \|\chi_{B_j}\|_{L^\varphi}^{-1}}{C_\gamma^{-1/\gamma}(\sum_{j=1}^k |\lambda_j|^\gamma)^{1/\gamma}}\right)\geq \varphi(B_{j_0}, \|\chi_{B_{j_0}}\|_{L^\varphi}^{-1})=1.$$
Otherwise, it follows from (\ref{the boundedness of sublinear operators on new Hardy spaces 1}) that
$$\sum_{j=1}^k \varphi\left(B_j, \frac{|\lambda_j| \|\chi_{B_j}\|_{L^\varphi}^{-1}}{C_\gamma^{-1/\gamma}(\sum_{j=1}^k |\lambda_j|^\gamma)^{1/\gamma}}\right)\geq \sum_{j=1}^k \frac{|\lambda_j|^\gamma}{\sum_{j=1}^k |\lambda_j|^\gamma}\varphi(B_j, \|\chi_{B_j}\|_{L^\varphi}^{-1})=1.$$

The above means that 
$$\Big(\sum_{j=1}^k |\lambda_j|^\gamma\Big)^{1/\gamma}\leq C_\gamma^{1/\gamma}\Lambda_q(\{\lambda_j a_j\}_{j=1}^k)\leq C\|f\|_{H^\varphi}.$$

Therefore, by assumption (i), we obtain that
$$\|Tf\|_{\mathcal B_\gamma}=\left\|T\Big(\sum_{j=1}^k \lambda_j a_j\Big)\right\|_{\mathcal B_\gamma}\leq C \Big(\sum_{j=1}^k |\lambda_j|^\gamma\Big)^{1/\gamma}\leq C\|f\|_{H^\varphi}.$$
Since $H^{\varphi,q,s}_{\rm fin}(\mathbb R^n)$ is dense in $H^\varphi(\mathbb R^n)$, a density argument gives the desired result.

The case (ii) is similar by using the fact that $H^{\varphi,\infty,s}_{\rm fin}(\mathbb R^n)\cap C(\mathbb R^n)$ is dense in $H^{\varphi,\infty,s}_{\rm fin}(\mathbb R^n)$ in the quasi-norm $\|\cdot\|_{H^\varphi}$, see the  below lemma.
\end{proof}

We end the paper by the following lemma.

\begin{Lemma}
Let $\varphi$ be a growth function satisfying uniformly locally dominated convergence condition, and $(\varphi,\infty,s)$ be an admissible triplet. Then, $H^{\varphi,\infty,s}_{\rm fin}(\mathbb R^n)\cap C^\infty(\mathbb R^n)$ is dense in $H^{\varphi,\infty,s}_{\rm fin}(\mathbb R^n)$ in the quasi-norm $\|\cdot\|_{H^\varphi}$.
\end{Lemma}

\begin{proof}
We take $q\in (q(\varphi),\infty)$ and $\phi\in\mathcal S(\mathbb R^n)$ satisfying supp $\phi\subset B(0,1)$, $\int_{\mathbb R^n}\phi(x)dx=1$. Then, the proof of the lemma is simple since it follows from the fact that for every $(\varphi,\infty,s)$-atom $a$ supported in ball $B(x_0,r)$,
$$\lim\limits_{t\to 0}\|a- a*\phi_t\|_{L^q_\varphi(B(x_0, 2r))}=0$$
as $\varphi$ satisfies uniformly locally dominated convergence condition.
\end{proof}

\subsection*{Acknowledgment}
The author would like to thank Prof. Aline Bonami, Prof. Sandrine Grellier and Prof. Dachun Yang for many helpful suggestions and discussions. He would also like to thank Prof. Sandrine Grellier for her carefully reading and revision of the manuscript. The author is deeply indebted to them.

\end{document}